\documentclass[12pt,leqno]{amsart}

\usepackage{amsmath,amssymb,amsthm}
\usepackage{hyperref}
\usepackage{a4wide}
\usepackage{color}

\newtheorem{theorem}{Theorem}[section]

\newtheorem{lemma}[theorem]{Lemma}
\newtheorem{corollary}[theorem]{Corollary}
\newtheorem{definition}[theorem]{Definition}
\theoremstyle{remark}
\newtheorem{remark}[theorem]{Remark}

\numberwithin{equation}{section}
\newcommand{\rr}{\mathbb{R}}
\newcommand{\zz}{\mathbb{Z}}
\newcommand{\Da}{\partial_{t}^{\alpha}}
\newcommand{\Dc}{D_{t}^{\alpha}}
\newcommand{\eps}{\varepsilon}

\newcommand{\sgn}{\operatorname{sgn}}
\newcommand{\BUC}{\mathrm{BUC}}
\newcommand{\EE}{\mathbb{E}}
\newcommand{\cD}{\mathcal{D}}
\makeatletter
\newcommand{\namedlabel}[2]{\phantomsection\def\@currentlabel{#2}\label{#1}}
\makeatother

\begin{document}

\title[Energy methods for time-fractional equations]{Energy methods and optimal decay for time-fractional equations}

\author[L. I. Ignat]{Liviu I. Ignat}
\address[L. I. Ignat]{Institute of Mathematics ``Simion Stoilow'' of the Romanian Academy, 21 Calea Grivitei Street, 010702 Bucharest, Romania.
\newline\indent
Academy of Romanian Scientists, Ilfov Street, no.~3, Bucharest, Romania.}
\email{liviu.ignat@gmail.com}


\begin{abstract}
We study energy methods for the time-fractional equation $\partial_t^\alpha(u-u_0)-(a(t,x)u_x)_x+(F(t,x,u))_x=Q(t,x,u)$ on the real line, where $0<\alpha<1$, the coefficient $a$ is measurable and bounded from above and from below by positive constants, $\sigma\partial_xF\ge0$, $\sigma Q\le0$, and the data are in $L^1\cap L^\infty$.
Using Yosida approximations of the kernel, an iteration on the $L^{2^m}$ norms gives the parabolic decay $\|u(t)\|_{L^\infty}\le C\|u_0\|_{L^1}t^{-\alpha/2}$.
 In $L^2$ it gives the rate $t^{-\alpha/4}$, which, already for the linear equation, improves the rate $t^{-\alpha/5}$ obtained before by the classical energy method.
For the convection--diffusion equation $\partial_t^\alpha(u-u_0)-u_{xx}+(|u|^{q-1}u)_x=0$, $q>1$, the parabolic rate is optimal for $q\ge2$, and for $1<q<2$ the convection improves it.
An $L^1$ time-regularity estimate of B\'enilan--Crandall type gives the pointwise flux bound $|u|^{q-1}u-u_x=O(\|u_0\|_{L^1}t^{-\alpha})$.
This bound takes the place of the Ole\u{\i}nik estimate and yields the optimal rate $\|u(t)\|_{L^\infty}\lesssim t^{-\alpha/q}$ for $1<q<2$. For nonnegative data, matching lower bounds hold for the $L^p$ norms.
\end{abstract}

\maketitle


\section{Introduction}

We study the decay of the solutions of the time-fractional equation
\begin{equation}\label{eq:general}
\Da(u-u_0)-\big(a(t,x)u_x\big)_x+\big(F(t,x,u)\big)_x=Q(t,x,u),\qquad t>0,\ x\in\rr,
\end{equation}
where $0<\alpha<1$ and $u_0$ is the initial datum. For $\beta>0$ we write
\[
g_\beta(t)=\frac{t^{\beta-1}}{\Gamma(\beta)},\qquad t>0.
\]
The operator $\Da$ is the Riemann--Liouville derivative of order $\alpha$,
\[
\Da v(t):=\frac{d}{dt}\int_0^tg_{1-\alpha}(t-s)\,v(s)\,ds=\frac{1}{\Gamma(1-\alpha)}\frac{d}{dt}\int_0^t\frac{v(s)}{(t-s)^{\alpha}}\,ds .
\]
Functions in $W^{1,1}(0,T)$, $T>0$, are identified with their continuous representatives. For such a function $v$, the Caputo derivative is
\[
\Dc v(t):=(g_{1-\alpha}\ast v')(t)=\frac{1}{\Gamma(1-\alpha)}\int_0^t\frac{v'(s)}{(t-s)^{\alpha}}\,ds=\Da\big(v-v(0)\big)(t),\quad\text{for a.e.\ }t\in(0,T) .
\]

The coefficient $a$ is measurable in $(t,x)$, with $\lambda\le a\le\Lambda$. The flux and the reaction term satisfy (the precise assumptions (\ref{hyp:A}), (\ref{hyp:F}) and (\ref{hyp:Q}) are stated in Section~\ref{sec:wp.general})
\[
\sigma\,\partial_xF(t,x,\sigma)\ge0,\qquad \sigma\,Q(t,x,\sigma)\le0 .
\]
 Throughout the paper we denote
\[
M:=\|u_0\|_{L^1(\rr)} .
\]
The model example that started this work is the convection--diffusion equation
\begin{equation}\label{convection}
\Da (u-u_0)-u_{xx}+(f(u))_x=0,\quad t>0,\ x\in\rr,
\qquad f(u)=|u|^{q-1}u,\quad q>1,
\end{equation}
which  {we analyze in more detail in Section~\ref{sec:L1}}.

Our main result concerns the weak solutions of \eqref{eq:general}, defined in Section~\ref{sec:wp.general}.

\begin{theorem}[Parabolic decay by energy methods]\label{thm:intro.energy}Assume \textup{(\ref{hyp:A})}, \textup{(\ref{hyp:F})} and \textup{(\ref{hyp:Q})}, let $u_0\in L^1\cap L^\infty(\rr)$ and let $u$ be the weak solution of \eqref{eq:general}. For a.e.\ $t>0$ and $1\le p\le\infty$,
\[
\|u(t)\|_{L^p(\rr)}\le C(\alpha)^{1-\frac1p}\lambda^{-\frac12(1-\frac1p)}M\,t^{-\frac\alpha2(1-\frac1p)}.
\]
\end{theorem}

Let us first explain what the main idea behind this result is. For simplicity we consider the linear equation
\begin{equation}\label{liniar}
\Da(u-u_0)-u_{xx}=0.
\end{equation}
It is  {well known} that its solutions satisfy (see \cite{Kemppainen2014})
\begin{equation}
  \label{eq:decay.liniar.1}
\|u(t)\|_{L^2(\rr)}\lesssim M\,t^{-\alpha/4},  
\end{equation}
and this rate is optimal.
The common source of difficulty in obtaining this decay by energy methods is the memory of $\Da$
{which remembers the initial value of the solution}.
The classical energy method for time-fractional equations (\cite{Kemppainen2014,KSZ17} in $\rr^d$, {\cite{MR3296607} on bounded domains}, \cite{DVV19,AV19} for abstract operators) and ours start in the same way. For a regular solution of \eqref{liniar} let
\[
N(t)=\|u(t)\|_{L^2(\rr)},\qquad D(t)=\|u_x(t)\|_{L^2(\rr)}.
\]
Multiplying the equation by $u$ and using the identity \eqref{acv2}, we obtain
\begin{equation}
  \label{eq:first.integration}
N\,\Da N+D^2\le g_{1-\alpha}(t)\int_\rr u(t)\,u_0\,dx.  
\end{equation}
The two methods differ only in the way they estimate the last term, which is the only one that contains the datum.
The classical method  
uses the Cauchy--Schwarz inequality and the identity
\[
\Da N-\|u_0\|_{L^2(\rr)}\,g_{1-\alpha}=\Da\big(N-\|u_0\|_{L^2(\rr)}\big),
\]
to obtain
\[
N\,\Da\big(N-\|u_0\|_{L^2(\rr)}\big)+D^2\le0.
\]
With Nash's inequality $N^6\le M^4D^2$ it becomes
\[
\Da\big(N-\|u_0\|_{L^2(\rr)}\big)+\frac{N^5}{M^4}\le0 .
\]
The comparison principle for  fractional ODEs and  Lemma~\ref{decay.odes}  give  \cite[Theorem 6.1]{Kemppainen2014},
\[
\|u(t)\|_{L^2(\rr)}\lesssim\|u_0\|_{L^2(\rr)}^{1/5}M^{4/5}\,t^{-\alpha/5}.
\]
 Our idea is to estimate the {right-hand} side in \eqref{eq:first.integration} without involving the $L^2$-norm of the initial data:
\begin{align*}
  N\,\Da N+D^2&\le g_{1-\alpha}(t)\,M\|u(t)\|_{L^\infty(\rr)}\leq g_{1-\alpha}(t)\,M N^{1/2}D^{1/2}\\
&\leq \frac12D^2+\big(g_{1-\alpha}(t)MN^{1/2}\big)^{4/3},
\end{align*}
and then,
\[
N\,\Da N+\frac{N^6}{2M^4}\le\big(M\,g_{1-\alpha}(t)\big)^{4/3}N^{2/3}.
\]
Only the mass of the datum appears in this inequality. Using well-constructed supersolutions gives the correct decay \eqref{eq:decay.liniar.1}.
This is the main idea of this paper that we apply to general nonlinear equations.

In one dimension the new energy method thus gives the rate $t^{-\alpha/4}$ in $L^2(\rr)$, also for coefficients that depend on time. For $a\equiv1$ this rate is optimal. For bounded coefficients, this answers the question raised in \cite{Kemppainen2014}, whether the rate $t^{-\alpha/5}$ of the energy method is optimal.

For \eqref{liniar}, the classical energy method can also be applied to each dyadic frequency block to obtain the optimal decay rate, see Section \ref{sec:linear}. 
This dyadic method can be adapted to the nonlinear equation \eqref{convection} for small initial data: if $q\ge2$ and $\|u_0\|_{L^1(\rr)}\|u_0\|_{L^\infty(\rr)}^{q-2}$ is small enough, it gives $\|u(t)\|_{L^2(\rr)}\lesssim\|u_0\|_{L^1(\rr)}t^{-\alpha/4}$.
The smallness enters because the convection does not vanish on a single block, and it couples the blocks. We will not pursue this issue in this paper.

 For particular equations,  the structure of the nonlinearity gives a faster decay, and \eqref{convection} is an example.  
The decay rate given by Theorem~\ref{thm:intro.energy} is optimal for $q\ge2$, but for $1<q<2$ the convection gives the faster hyperbolic rate $t^{-\alpha/q}$.
In this case we first get that   $A(u)=u_{xx}-f(u)_x$ satisfies
\[
\Gamma(1-\alpha)\|A(u(t))\|_{L^1(\rr)}\le\frac{\|u(t)-u_0\|_{L^1(\rr)}}{t^\alpha}+\alpha\int_0^t\frac{\|u(t)-u(s)\|_{L^1(\rr)}}{(t-s)^{1+\alpha}}\,ds .
\]
The first term is at most $2Mt^{-\alpha}$. The second one requires a bound on the time increments of $u$ in $L^1$.
This is obtained in Theorem~\ref{thm:intro.L1}, and hence $\|A(u(t))\|_{L^1(\rr)}\lesssim Mt^{-\alpha}$.
This turns into a pointwise bound for the flux $|u|^{q-1}u-u_x$, which finally gives the correct decay of the solutions.

In Theorems~\ref{thm:intro.L1} and \ref{thm:intro.decay}, $u$ is the mild solution of \eqref{convection} constructed in Section~\ref{sec:wp.mild}, and $\omega_\alpha$ is defined by
\[
\omega_\alpha=1+\frac{\omega_\alpha^\alpha}{1-\alpha},\qquad \omega_\alpha>1.
\]

\begin{theorem}[$L^1$ time regularity and flux bound]\label{thm:intro.L1}
Let $q>1$ and $u_0\in L^1\cap L^\infty(\rr)$. For all $0\le s<\tau$ and $t>0$,
\[
\|u(\tau)-u(s)\|_{L^1(\rr)}\le 2\omega_\alpha M\Big(1-\big(\tfrac s\tau\big)^\alpha\Big),\]
\[
\big\|\big(|u|^{q-1}u-u_x\big)(t)\big\|_{L^\infty(\rr)}\le \frac{\omega_\alpha M}{\Gamma(1-\alpha)}\,t^{-\alpha}.
\]
\end{theorem}

\begin{theorem} \label{thm:intro.decay}
Let $q>1$ and $u_0\in L^1\cap L^\infty(\rr)$. For all $t>0$,
\[
\|u(t)\|_{L^\infty(\rr)}\le\min\Big\{\big(\omega_\alpha Mg_{1-\alpha}(t)\big)^{1/q},\ M\big(2\omega_\alpha g_{1-\alpha}(t)\big)^{1/2}\Big\}.
\]
If $u_0\ge0$ and $u_0\neq0$, then for $t\ge1$ and $1<p\le\infty$,
\[
\|u(t)\|_{L^p(\rr)}\simeq t^{-\frac{\alpha}{\min\{q,2\}}(1-\frac1p)},
\]
with constants depending only on $\alpha$, $q$ and $u_0$.
\end{theorem}

Theorem~\ref{thm:intro.L1} is an estimate of B\'enilan--Crandall type \cite{BC81} (see also \cite{Hauer21} and the references therein), but the operator is not homogeneous, and the convection is not a Lipschitz perturbation in $L^1(\rr)$ as in \cite{Hauer21}. For the time-fractional porous medium equation $\Dc u=\Delta(u^m)$, $m\neq1$, on bounded domains of $\rr^N$, $N\ge3$, with zero Dirichlet conditions, the $L^p$--$L^\infty$ smoothing effect $\|u(t)\|_{L^\infty}\le\kappa\|u_0\|_{L^p}^{1/m}t^{-\alpha/m}$, for $p\ge1+m$ and $p>N/2$, is proved in \cite[Theorem 1.4]{BGI26}. The proof uses a time monotonicity estimate of B\'enilan--Crandall type, based on the homogeneity of the operator, which is not available here.

In the case $\alpha=1$, that is, for
\[
u_t-u_{xx}+(|u|^{q-1}u)_x=0,
\]
the theory is classical. The $L^1$ norm is nonincreasing, the convective term does not contribute to the energy identity, and
\[
\|u(t)\|_{L^\infty(\rr)}\lesssim t^{-1/2}\qquad\text{for all }q>1\text{ and all data in }L^1\cap L^\infty(\rr).
\]
Estimates of this type were obtained by Schonbek with the Fourier splitting method \cite{Schonbek1980}, and in this generality by Escobedo and Zuazua \cite{EZ} with energy methods (see also \cite[Ch.~5, Theorem~5.2, p.~59]{Zuazua25}).
The {large-time} behavior depends on $q$ \cite{EZ,EVZ}. For $q>2$ the solution behaves like the heat kernel,
for $q=2$ like a self-similar solution of the Burgers equation, and for $1<q<2$ the hyperbolic part dominates.
For $1<q<2$ the optimal $L^\infty$ decay rests on the one-sided Ole\u{\i}nik estimate \cite{Ole63}, see \cite[Lemmas~1.1 and 1.2, p.~48]{EVZ}. For $q>2$, a one-sided H\"older estimate for nonnegative solutions was established by De Nitti, Ignat and Quir\'os \cite{DIQ}.
For $0<\alpha<1$ these one-sided estimates are replaced by the two-sided flux bound of Theorem~\ref{thm:intro.L1}, which depends on the data only through the mass.

The multidimensional analogues of these results, as well as problems with fractional diffusion in space, will be analyzed elsewhere. Since including them would make the paper more technical, we prefer to present here the one-dimensional case, which captures the main idea behind the energy methods for fractional-in-time evolution problems. 

The paper is organized as follows. Section~\ref{sec:prelim} collects tools from fractional calculus, including the Yosida approximations of the kernel.
Section~\ref{sec:wp} states the well-posedness theory, for mild and regular solutions of \eqref{convection} and for weak solutions of \eqref{eq:general}, and it proves the $L^1$ theory of regular solutions. The theory of Volterra equations with m-accretive operators gives generalized solutions and the $L^1$ contraction \cite{CoGL96,Gripenberg85}, but we did not find a reference that gives the regularity of Theorem~\ref{thm:wp} or treats the nonlinear operators of Theorem~\ref{thm:wp.general}. This is why we decided to 
  sketch the proof of Theorem~\ref{thm:wp} in Appendix~\ref{app:wp} and to prove Theorem~\ref{thm:wp.general} in Appendix~\ref{app:general}.
Section~\ref{sec:linear} recalls the classical energy method for the linear equation and proves the optimal $L^2$ decay with its dyadic version.
Section~\ref{sec:energy.nonlinear} proves Theorem~\ref{thm:intro.energy}, the key estimate being the comparison with barriers of Lemma~\ref{lem:abstract}.
Section~\ref{sec:L1} proves Theorems~\ref{thm:intro.L1} and \ref{thm:intro.decay}, which for \eqref{convection} with $1<q<2$ give a faster decay than Theorem~\ref{thm:intro.energy}.

\section{Preliminaries}\label{sec:prelim}

In this section we collect the tools from fractional calculus used in the rest of the paper: an identity for the Caputo derivative, an energy inequality, comparison results for fractional ODEs, and the Yosida approximations of the kernel $g_{1-\alpha}$.

\subsection{Fractional derivatives}
In what follows we will denote
\[
(f\ast g)(t)=\int_0^t f(t-s)g(s)\,ds.\]
For functions on $\rr$, we denote by $(f\star g)(x)=\int_\rr f(x-y)g(y)\,dy$ the convolution in $x$.
We use the identity of \cite{Allen2015} in the following form. If $v\in W^{1,1}(0,t)$ is Lipschitz on $[t/2,t]$, then
\begin{align}\label{acv2}
\Gamma(1-\alpha)\Da\big(v-v(0)\big)(t)&=\Gamma(1-\alpha)\,(g_{1-\alpha}\ast v')(t)\\
&=\frac{v(t)-v(0)}{t^\alpha}+\alpha\int_0^t\frac{v(t)-v(s)}{(t-s)^{1+\alpha}}\,ds,\notag
\end{align}
where the integral converges absolutely. Here $\Da$ at the point $t$ is a left derivative, since $g_{1-\alpha}\ast v'$ is continuous on $[3t/4,t]$.

\subsection{Fundamental inequalities}
Let $\mathcal H$ be a real Hilbert space, $T>0$, and let $v\in C([0,T];\mathcal H)$ satisfy
\[
v(0)=v_0,\qquad \partial_tv\in C((0,T];\mathcal H),\qquad \|\partial_tv(t)\|\le Ct^{\alpha-1},
\]
so that $\Da(v-v_0)=g_{1-\alpha}\ast\partial_tv$. Then for a.e.\ $t\in(0,T)$
\begin{equation}\label{eq:I}
\langle v,\Da(v-v_0)\rangle \ge \|v\|\,\Da\big(\|v\|-\|v_0\|\big).
\end{equation}
For nonnegative and nonincreasing kernels $k\in W^{1,1}(0,T)$ in place of $g_{1-\alpha}$, inequality \eqref{eq:I} is \cite[Lemma~3.3, p.~222]{MR3296607} (see \cite[Lemma 6.2]{Kemppainen2014} for $\mathcal H=L^2(\Omega)$).
For the singular kernel and the functions above, it follows from \eqref{acv2} in $\mathcal H$ and the first of the elementary inequalities
\begin{equation}\label{eq:elem}
\langle y,y-z\rangle\ge\|y\|\big(\|y\|-\|z\|\big),
\qquad \langle y,y-z\rangle\ge\tfrac12\big(\|y\|^2-\|z\|^2\big),
\qquad y,z\in\mathcal H .
\end{equation}
The second one is used in Appendix~\ref{app:general}.
Inequality \eqref{eq:I} applies to the regular solutions of Section~\ref{sec:wp}, by (\ref{R2}), and to the dyadic blocks of the solutions of the linear equation (Section~\ref{sec:linear}).

In one dimension we also use Agmon's inequality and the Nash inequality that it implies,
\begin{align}
\label{agmon}
\|v\|_{L^\infty(\rr)}^2&\le\|v\|_{L^2(\rr)}\|v_x\|_{L^2(\rr)}, & &v\in H^1(\rr),\\
\label{nash}
\|v\|^6_{L^2(\rr)}&\le\|v\|^4_{L^1(\rr)}\|v_x\|^2_{L^2(\rr)}, & &v\in L^1\cap H^1(\rr).
\end{align}

\subsection{Mittag-Leffler functions and comparison}
For $\beta>0$ and $b\in\rr$ the Mittag-Leffler functions are
\[
E_{\beta,b}(z)=\sum_{k\ge0}\frac{z^k}{\Gamma(\beta k+b)},\qquad E_\beta=E_{\beta,1},\qquad \frac1{\Gamma(-n)}=0\quad(n\in\mathbb N_0).
\]
For $0<\alpha<1$ and $x\ge 0$,
\begin{equation}\label{eq:ML}
0<E_\alpha(-x)\le \frac{1}{1+x/\Gamma(1+\alpha)}\le\frac{1}{1+x},\qquad
E_{\alpha,\alpha}(-x)\ge0.
\end{equation}
The bound $E_\alpha(-x)\le(1+x/\Gamma(1+\alpha))^{-1}$ is \cite[Theorem~4]{Simon}; see also \cite[Lemma~6.1 and (39), pp.~226--227]{MR3296607} and \cite[Ch.~3, Theorem~3.6, p.~71]{Jin}. For completeness we recall its short proof.
The function $\varphi(t)=E_\alpha(-t^\alpha)$ is nonincreasing, since $E_\alpha(-x)$ is completely monotone \cite[Ch.~3, Theorem~3.5, p.~69]{Jin}, and $\varphi+g_\alpha\ast\varphi=1$ \cite[Ch.~3, (3.22), p.~71]{Jin}. Hence
\[
1=\varphi(t)+\int_0^tg_\alpha(t-s)\varphi(s)\,ds\ge\varphi(t)\Big(1+\frac{t^\alpha}{\Gamma(1+\alpha)}\Big),
\]
which is this bound with $x=t^\alpha$. We recall that $\Gamma(1+\alpha)\le1$.
The positivity of
\[
E_{\alpha,\alpha}(-x)=-\alpha\frac{d}{dx}E_\alpha(-x)
\]
\cite[Ch.~3, (3.21), p.~71]{Jin} follows from the complete monotonicity of $E_\alpha(-x)$ \cite[Ch.~3, Theorem~3.5 and Corollary~3.3, pp.~69 and 71]{Jin}.
For $\nu\ge0$ the resolvent kernel of $\Da+\nu$ is
\[
r_\nu(t):=t^{\alpha-1}E_{\alpha,\alpha}(-\nu t^\alpha),\qquad t>0 .
\]
By \eqref{eq:ML}, $r_\nu\ge0$, and $r_0=g_\alpha$.

\begin{lemma}[Linear comparison]\label{lem:comparison}
Let $\nu\ge 0$, $w_0\in\rr$ and $h\in L^1_{loc}([0,\infty))$, and let $w\in L^\infty_{loc}([0,\infty))$, with $g_{1-\alpha}\ast(w-w_0)\in W^{1,1}_{loc}([0,\infty))$, satisfy
\[
\Da(w-w_0)+\nu w\le h\qquad\text{a.e.\ in }(0,\infty).
\]
Then
\[
w(t)\le w_0E_\alpha(-\nu t^\alpha)+\int_0^t r_\nu(t-s)h(s)\,ds,\qquad \text{for a.e.\ }t>0.
\]
\end{lemma}

\begin{proof}
Denote by $\mathcal L$ the Laplace transform. Since
\[
\mathcal L\big[E_\alpha(-\nu t^\alpha)\big](z)=\frac{z^{\alpha-1}}{z^\alpha+\nu},\qquad \mathcal L[r_\nu](z)=\frac{1}{z^\alpha+\nu},
\]
the function $\bar w=w_0E_\alpha(-\nu t^\alpha)+r_\nu\ast h$ solves $\Da(\bar w-w_0)+\nu \bar w=h$.
The difference $v=w-\bar w$ satisfies
\[
\partial_t(g_{1-\alpha}\ast v)+\nu v=-\eta,\qquad \eta\ge0 .
\]
Since $g_{1-\alpha}\ast v\in W^{1,1}_{loc}([0,\infty))$ vanishes at $t=0$, convolving with $g_\alpha$ gives $v+\nu g_\alpha\ast v=-g_\alpha\ast\eta$, whence $v=-r_\nu\ast\eta\le0$.
\end{proof}

The following result is \cite[Theorem 7.1]{MR3296607}, with the dependence on the parameters made explicit.

\begin{lemma}\label{decay.odes}
Let $\alpha\in (0,1)$, $\gamma,\nu>0$ and $v_{0}>0$, and let $v\in W^{1,1}_{loc}([0,\infty))$, with $v(0)=v_0$, be the nonnegative solution of
\begin{equation}\label{non.ode}
\Da (v-v_{0})+\nu v^{\gamma}=0, \quad t\ge 0.
\end{equation}
There exist constants $0<c_1\le c_2$, depending only on $\alpha$ and $\gamma$, such that
\[
\frac{c_1v_0}{1+(\nu v_0^{\gamma-1}t^\alpha)^{1/\gamma}}\le v(t)\le\frac{c_2v_0}{1+(\nu v_0^{\gamma-1}t^\alpha)^{1/\gamma}},\qquad t\ge0 .
\]
\end{lemma}

\begin{proof}
Since the derivative is {linear and} homogeneous under time dilations,
 it is sufficient to consider $\nu=1$. The sub- and supersolutions constructed in \cite[Section 7]{MR3296607} give
the result in this case 
and scaling back gives the claim.
\end{proof}

Lemma~\ref{decay.odes} exhibits a memory effect: $\nu v(t)^\gamma\approx v_0\,g_{1-\alpha}(t)$ for large $t$. This explains the rates of the classical energy method in the introduction.

\subsection{Yosida approximations of the kernel}\label{sec:yosida}
Let $X$ be a Banach space, $k\in W^{1,1}(0,T)$ and $N\in L^1(0,T;X)$. Then $\partial_t(k\ast N)=k(0)N+k'\ast N$, and since $k(0)=k(t)-\int_0^tk'(t-s)\,ds$,
\begin{equation}\label{eq:Rk}
\partial_t(k\ast N)(t)=k(t)N(t)+\int_0^t\big(N(t)-N(s)\big)\big(-k'(t-s)\big)\,ds=:\cD_k[N](t)\quad\text{for a.e.\ }t .
\end{equation}
Unlike \eqref{acv2}, this identity requires no regularity of $N$ in time. If $N$ is a bounded measurable function defined at every $t$, the right-hand side of \eqref{eq:Rk} is defined for every $t$, and we use it as the definition of $\cD_k[N](t)$.

The Yosida approximations of $\partial_t(g_{1-\alpha}\ast\cdot)$ are the operators $\partial_t(k_n\ast\cdot)$, with the kernels $k_n$ of the following lemma \cite[Section 2]{Zacher09}. We use them for the weak solutions of \eqref{eq:general}, which need not be regular in time.

\begin{lemma}\label{lem:yosida}
For $n\in\mathbb N$ let
\[
s_n(t)=E_\alpha(-nt^\alpha),\qquad k_n=ns_n,\qquad \Pi_\alpha=\Gamma(1+\alpha)\Gamma(1-\alpha)=\frac{\pi\alpha}{\sin(\pi\alpha)}\ge1 .
\]
Let $T>0$ and let $X$ be a Banach space.
\begin{enumerate}
\item[(a)] $k_n$ is positive and nonincreasing, $k_n(0)=n$, $k_n\in W^{1,1}(0,T)$, and $k_n\to g_{1-\alpha}$ in $L^1(0,T)$. Moreover, for $t>0$,
\[
k_n(t)\le \Pi_\alpha\,g_{1-\alpha}(t),\qquad -k_n'(t)\le4\Pi_\alpha\,|g_{1-\alpha}'(t)| .
\]
\item[(b)] $g_\alpha\ast k_n=1-s_n$. If $w\in L^1(0,T;X)$ and $\psi=\partial_t(k_n\ast w)$, then $w=n^{-1}\psi+g_\alpha\ast\psi$.
\end{enumerate}
\end{lemma}

\begin{proof}
The Laplace transform $\mathcal L[s_n](z)=z^{\alpha-1}/(z^\alpha+n)$ gives
\[
s_n+ng_\alpha\ast s_n=1,\qquad\text{that is,}\qquad g_\alpha\ast k_n=1-s_n.
\]
Since $g_{1-\alpha}\ast g_\alpha=1$, these are the kernels of \cite[Section 2]{Zacher09} for the pair $(g_{1-\alpha},g_\alpha)$. Hence $\partial_t(k_n\ast\cdot)$ are the Yosida approximations of $\partial_t(g_{1-\alpha}\ast\cdot)$, $k_n\in W^{1,1}(0,T)$, and $k_n\to g_{1-\alpha}$ in $L^1(0,T)$ \cite[Section 2]{Zacher09}.

(a) Let $e(x)=E_\alpha(-x)$, which is completely monotone \cite[Ch.~3, Theorem~3.5, p.~69]{Jin}. Hence $s_n$ is positive and nonincreasing, and $-e'$ is nonincreasing.
With $x=nt^\alpha$, \eqref{eq:ML} gives
\[
t^{\alpha}k_n(t)=x\,e(x)\le\Gamma(1+\alpha).
\]
Since $\frac x2(-e'(x))\le e(\frac x2)-e(x)\le e(\frac x2)$, the same bound gives
\[
-t^{1+\alpha}k_n'(t)=\alpha x^2\big(-e'(x)\big)\le2\alpha  x\,e\big(\tfrac x2\big)\le4\alpha\Gamma(1+\alpha) 
.
\]
(b) Since $k_n\ast w\in W^{1,1}(0,T;X)$ has value $0$ at $t=0$, we have \[
(k_n\ast w)(t)=\int_0^t \psi(s)ds=1\ast\psi.\]
 Hence
\[
1\ast\big(n^{-1}\psi+g_\alpha\ast\psi\big)=n^{-1}k_n\ast w+g_\alpha\ast k_n\ast w=s_n\ast w+(1-s_n)\ast w=1\ast w,
\]
and $w=n^{-1}\psi+g_\alpha\ast\psi$.
\end{proof}

\section{Well-posedness}\label{sec:wp}

This section states the existence theory used in the rest of the paper, first for \eqref{convection} and then for \eqref{eq:general}.
 The proof of Theorem~\ref{thm:wp} follows standard arguments for semilinear subdiffusion equations, and we sketch it in Appendix~\ref{app:wp}. We prove Theorem~\ref{thm:wp.general} in Appendix~\ref{app:general}.

For \eqref{convection} we call \emph{flux} a function $f\in C^1(\rr)$ with $f(0)=0$, and \emph{regular flux} a flux of class $C^2$ with bounded derivative.
We use two notions of solution. For regular fluxes and smooth data we use \emph{regular} solutions, for which all the computations of Sections~\ref{sec:linear} and \ref{sec:L1} are justified.
For the power flux and data in $L^1\cap L^\infty(\rr)$ the solution is the \emph{mild} solution, which is the limit of regular solutions with approximating fluxes and data.
Our estimates for regular solutions hold for every regular flux, and they depend on the data only through norms $\|u_0\|_{L^p(\rr)}$. Hence they pass to the limit (Remark~\ref{rem:transfer}).

\subsection{The kernels of the linear equation}\label{sec:kernels}
We use the Fourier transform
\[
\hat g(\xi)=(2\pi)^{-1/2}\int_\rr e^{-ix\xi}g(x)\,dx.
\]
For $t>0$ let $S(t),P(t)\in L^1(\rr)$ be the functions such that
\[
\widehat{S(t)\star g}(\xi)=E_\alpha(-\xi^2t^\alpha)\,\hat g(\xi),\qquad \widehat{P(t)\star g}(\xi)=t^{\alpha-1}E_{\alpha,\alpha}(-\xi^2t^\alpha)\,\hat g(\xi).
\]
They are the fundamental solutions $G_\alpha$ and $\mathcal G_\alpha$ of \cite[Ch.~7, Proposition~7.1 and its proof, p.~278]{Jin}, with $d=1$.
We use the following properties of these kernels.
\begin{itemize}
\item[(K1)]\namedlabel{K1}{K1} $S(t)$ and $P(t)$ are positive, with
\[
\int_\rr S(t,x)\,dx=1,\qquad\int_\rr P(t,x)\,dx=g_\alpha(t)
\]
\cite[Ch.~7, Lemma~7.1, p.~279]{Jin}. Hence, for $1\le p\le\infty$,
\[
\|S(t)\star g\|_{L^p(\rr)}\le\|g\|_{L^p(\rr)},\qquad \|P(t)\star g\|_{L^p(\rr)}\le g_\alpha(t)\|g\|_{L^p(\rr)} .
\]
\item[(K2)]\namedlabel{K2}{K2} $P(t)\in W^{1,1}(\rr)$, with
\[
\|\partial_xP(t)\|_{L^1(\rr)}=g_{\alpha/2}(t)
\]
\cite[Ch.~7, proof of Lemma~7.3, p.~282]{Jin}. The maps $t\mapsto S(t),P(t),\partial_xP(t)$ are continuous from $(0,\infty)$ into $L^1(\rr)$, and $S(t)\star g\to g$ in $L^1(\rr)$ as $t\to0$ for every $g\in L^1(\rr)$.
\end{itemize}
The last two properties follow from the continuity of dilations in $L^1(\rr)$ and from the scaling
\begin{equation}\label{eq:scaling}
S(t,x)=t^{-\alpha/2}S(1,t^{-\alpha/2}x),\qquad P(t,x)=t^{\alpha/2-1}P(1,t^{-\alpha/2}x),
\end{equation}
which follows from the definition of $S(t)$ and $P(t)$.

\subsection{Mild and regular solutions of \eqref{convection}}\label{sec:wp.mild}
A \emph{mild solution} on $[0,T]$ with flux $f$ and datum $u_0\in L^1\cap L^\infty(\rr)$ is a function $u\in C([0,T];L^1(\rr))$ with $\sup_{t\le T}\|u(t)\|_{L^\infty(\rr)}<\infty$ such that
\begin{equation}\label{eq:mild}
u(t)=S(t)\star u_0-\int_0^t\partial_xP(t-s)\star f(u(s))\,ds,\qquad 0\le t\le T .
\end{equation}
Here $S(0)\star g=g$, and the integral is a Bochner integral in $L^1(\rr)$, which converges by (\ref{K2}) since $f(u)\in C([0,T];L^1(\rr))$.
A \emph{global mild solution} is a mild solution on $[0,T]$ for every $T>0$.

\begin{definition}[Regular solutions]\label{def:regular}
A \emph{regular solution} on $[0,T]$ with a regular flux $f$ is a function $u$ with the following properties. Its datum is $u_0:=u(0)$.
\begin{enumerate}
\item[(R1)]\namedlabel{R1}{R1} $u\in C([0,T];W^{2,1}(\rr))$;
\item[(R2)]\namedlabel{R2}{R2} $u\in C^1((0,T];L^1\cap L^\infty(\rr))$, the partial derivative $\partial_tu$ exists and is continuous on $(0,T]\times\rr$, and
\[
\|\partial_tu(t)\|_{L^1(\rr)}+\|\partial_tu(t)\|_{L^\infty(\rr)}\le Ct^{\alpha-1},\qquad 0<t\le T;
\]
\item[(R3)]\namedlabel{R3}{R3} the equation holds in the two forms
\[
\big(g_{1-\alpha}\ast\partial_tu(\cdot,x)\big)(t)=u_{xx}(t,x)-(f(u))_x(t,x),\qquad u(t)=u_0+\big(g_\alpha\ast(u_{xx}-f(u)_x)\big)(t),
\]
the first for every $t\in(0,T]$ and a.e.\ $x\in\rr$, the second in $L^1(\rr)$ for every $t\in[0,T]$.
\end{enumerate}
\end{definition}

Condition (\ref{R1})   justifies the integrations by parts in $x$ and the evaluation of $u$ and $u_x$ at a point.
The pointwise bound in (\ref{R2}) is used in two ways. First, it allows us to apply \eqref{acv2} at every fixed $x$ (Lemma~\ref{lem:kato} and Step 1 of the proof of Theorem~\ref{thm:L1}): by (\ref{R1})--(\ref{R2}), every $u(\cdot,x)$ belongs to $W^{1,1}(0,T)$ and is Lipschitz away from $t=0$.
Second, together with the $L^1$ bound it gives
\[
\partial_tu\in C((0,T];L^2(\rr)),\qquad \|\partial_tu(t)\|_{L^2(\rr)}\le Ct^{\alpha-1},
\]
as the energy inequality \eqref{eq:I} requires.

Regular solutions require smooth data and a regular flux. The term $f''(u)u_x^2$ appears when the equation is differentiated twice in $x$, and the bound on $f'$ gives global solutions.
The power flux $|\sigma|^{q-1}\sigma$ is not regular: it is of class $C^2$ only for $q>2$, and its derivative is unbounded. We approximate it as follows.
Let $f$ be a flux and $R>0$. A sequence $(f_n)_{n\ge1}$ of regular fluxes is an \emph{admissible approximation of $f$ on $[-R,R]$} if
\[
\delta_n(R):=\max_{|\sigma|\le R}|f_n'(\sigma)-f'(\sigma)|\to0 .
\]
Then
\[
L_f(R):=\sup_n\max_{|\sigma|\le R}|f_n'(\sigma)|<\infty,
\]
a constant that depends on the sequence $(f_n)$, and, for $|\sigma|,|\sigma'|\le R$,
\begin{equation}\label{eq:adm}
|f_n(\sigma)-f(\sigma)|\le\delta_n(R)|\sigma|,\qquad |f_n(\sigma)-f_n(\sigma')|\le L_f(R)|\sigma-\sigma'| .
\end{equation}
For instance, let $\varrho_n=n\varrho(n\,\cdot)$, with $\varrho\in C_c^\infty(-1,1)$ a nonnegative mollifier of integral one, and let $\chi\in C_c^\infty(\rr)$ equal $1$ on $[-R-1,R+1]$. Then
\[
f_n(\sigma)=\int_0^\sigma\big(\varrho_n\star (\chi f')\big)(r)\,dr
\]
is an admissible approximation of $f$ on $[-R,R]$.

\begin{theorem}[Well-posedness]\label{thm:wp}The following well-posedness results hold.
\begin{enumerate}
\item[(i)] \emph{(Regular solutions.)} Let $f$ be a regular flux and $u_0\in W^{3,1}(\rr)$ (e.g.\ $u_0\in C_c^\infty(\rr)$).
There is a unique global mild solution, and it is a regular solution on every $[0,T]$.
\item[(ii)] \emph{(Mild solutions.)} Let $f$ be a flux (in particular $f(u)=|u|^{q-1}u$), and $u_0\in L^1\cap L^\infty(\rr)$.
There is a unique global mild solution $u\in C([0,\infty);L^1(\rr))$ which satisfies for all $t\geq 0$
\[
\begin{gathered}
 \int_\rr u(t)\,dx=\int_\rr u_0\,dx,\\
 -\|(u_0)_-\|_{L^\infty(\rr)}\le u(t,x)\le\|(u_0)_+\|_{L^\infty(\rr)},\quad \text{a.e.\ }x\in \rr,\\
 \|u(t)\|_{L^p(\rr)}\le\|u_0\|_{L^p(\rr)},\quad 1\le p\le\infty.
\end{gathered}
\]
Moreover, if $v$ is the mild solution with datum $v_0$, then for $t\ge0$,
\[
\|u(t)-v(t)\|_{L^1(\rr)}\le\|u_0-v_0\|_{L^1(\rr)}.
\]
\item[(iii)] \emph{(Approximation.)} In (ii), let $R\ge\|u_0\|_{L^\infty(\rr)}$ be positive, let $(f_n)$ be an admissible approximation of $f$ on $[-R,R]$, and let $v_{0,n}\in C_c^\infty(\rr)$ satisfy
\[
\|v_{0,n}\|_{L^\infty(\rr)}\le R,\qquad v_{0,n}\to u_0\quad\text{in }L^1(\rr).
\]
Then the regular solutions $v_n$ with flux $f_n$ and datum $v_{0,n}$ satisfy $|v_n|\le R$ and
\[
\sup_{0\le t\le T}\|v_n(t)-u(t)\|_{L^p(\rr)}\to0,\qquad 1\le p<\infty,\ T>0.
\]
 \end{enumerate}
\end{theorem}

\begin{remark}\label{stantard.data}
 Truncation and convolution give an approximating sequence  $u_{0,n}=\varrho_n\star (1_{[-n,n]}u_0)$ that satisfies
\[
u_{0,n}\in C_c^\infty(\rr),\qquad u_{0,n}\to u_0\ \text{in }L^1(\rr),\qquad \|u_{0,n}\|_{L^p(\rr)}\le\|u_0\|_{L^p(\rr)},\quad 1\le p\le\infty.
\]
Hence they are admissible in (iii) for every $R\ge\|u_0\|_{L^\infty(\rr)}$, and they are nonnegative if $u_0$ is.
\end{remark}

\subsection{The \texorpdfstring{$L^1$}{L1} theory for regular solutions}\label{sec:L1theory}
In this subsection we prove the $L^1$ properties of regular solutions and explain how estimates pass to mild solutions. We write
\[
A(u):=u_{xx}-f(u)_x,\qquad \sgn_+(r)=\begin{cases}1,&r>0,\\0,&r\le0,\end{cases}
\]
so that (\ref{R3}) reads $g_{1-\alpha}\ast\partial_tu=A(u)$, and we set $\sgn(0)=0$.
The following lemma is a pointwise-in-time form, for the kernel $g_{1-\alpha}$, of the inequality for the positive part of a subsolution \cite[Lemma 4.1]{MR3296607} and of the convexity inequality \cite[Lemma 2.1]{WWZ21}.

\begin{lemma}[Kato inequality for the Caputo derivative]\label{lem:kato}
Let $\phi_1$ have the regularity (\ref{R1})--(\ref{R2}) on $[0,T]$, let $\ell\ge0$, and let
\[
\phi=\phi_1-\ell,\qquad \Phi_+(t)=\int_\rr\phi(t,x)_+\,dx .
\]
Then $\Phi_+\in W^{1,1}(0,T)$ is locally Lipschitz on $(0,T]$, and
\[
\int_\rr\sgn_+(\phi(t,x))\,\big(g_{1-\alpha}\ast\partial_t\phi(\cdot,x)\big)(t)\,dx\ \ge\ (g_{1-\alpha}\ast\Phi_+')(t),\qquad 0<t\le T.
\]
If $\ell=0$, the same holds with $\sgn$ in place of $\sgn_+$ and $\Phi(t)=\|\phi(t)\|_{L^1(\rr)}$ in place of $\Phi_+$.
\end{lemma}

\begin{proof}
Since $\phi_+\le|\phi_1|$, the function $\Phi_+$ is finite, and by (\ref{R2})
\[
|\Phi_+(t)-\Phi_+(s)|\le\|\phi_1(t)-\phi_1(s)\|_{L^1(\rr)}\le\int_s^t\|\partial_t\phi_1(r)\|_{L^1(\rr)}\,dr\le \frac C\alpha\,(t^\alpha-s^\alpha),
\]
which gives the regularity of $\Phi_+$.
  We multiply \eqref{acv2} by $\sgn_+(\phi(t,x))$ and use the inequality $\sgn_+(r)(r-r')\ge r_+-r'_+$:
\[
\Gamma(1-\alpha)\sgn_+(\phi(t,x))(g_{1-\alpha}\ast\partial_t\phi)(t,x)\ge\frac{\phi(t,x)_+-\phi(0,x)_+}{t^\alpha}+\alpha\int_0^t\frac{\phi(t,x)_+-\phi(s,x)_+}{(t-s)^{1+\alpha}}\,ds .
\]
We integrate in $x$. Fubini's theorem applies because
\[
\int_0^t\frac{\|\phi_1(t)-\phi_1(s)\|_{L^1(\rr)}}{(t-s)^{1+\alpha}}\,ds<\infty,
\]
and the right-hand side becomes the right-hand side of \eqref{acv2} for $v=\Phi_+$.
For $\ell=0$ and $\sgn$ we use $\sgn(r)(r-r')\ge|r|-|r'|$ instead.
\end{proof}

\begin{lemma}[$L^1$ dissipativity]\label{lem:dissip}
Let $f$ be a flux and $u\in W^{2,1}(\rr)$, and let $v$ be either in $W^{2,1}(\rr)$ or a nonnegative constant $\ell$ (then $A(v)=0$). Then
\[
\int_\rr\sgn_+(u-v)\,\big(A(u)-A(v)\big)\,dx\le0 .
\]
If $v\in W^{2,1}(\rr)$, the same holds with $\sgn$ in place of $\sgn_+$.
\end{lemma}

\begin{proof}
This is the computation in the proof of the $L^1$ contraction for viscous conservation laws \cite[Ch.~VI, proof of Theorem~6.3.2, pp.~218--219]{Dafermos}. Let $\phi=u-v$, and let $\eta_\eps(r)=0$ for $r\le0$, $\eta_\eps(r)=r^2/(4\eps)$ for $0<r\le2\eps$ and $\eta_\eps(r)=r-\eps$ for $r>2\eps$. Integrating by parts (if $v$ is a positive constant, $\eta_\eps'(\phi)$ has compact support because $u$ vanishes at infinity),
\[
\int_\rr\eta_\eps'(\phi)\big(A(u)-A(v)\big)\,dx=-\int_\rr\eta_\eps''(\phi)\phi_x^2\,dx+\int_\rr\eta_\eps''(\phi)\big(f(u)-f(v)\big)\phi_x\,dx .
\]
The first term on the right is nonpositive. If $L$ is the Lipschitz constant of $f$ on the convex hull of the ranges of $u$ and $v$, then $|\eta_\eps''(\phi)(f(u)-f(v))|\le L$, this function tends to $0$ pointwise, and $\phi_x\in L^1(\rr)$. Dominated convergence on both sides, with $\eta_\eps'(\phi)\to\sgn_+(\phi)$, gives the claim.
Adding the inequalities for $(u,v)$ and for $(v,u)$ gives the second one.
\end{proof}

\begin{lemma}[Comparison]\label{lem:L3}
Let $\Phi\in W^{1,1}(0,T)$ with $\Phi(0)=0$, and $h\in L^1(0,T)$. Then
\[
g_{1-\alpha}\ast\Phi'\le h\ \text{ a.e.\ on }(0,T)\quad\Longrightarrow\quad\Phi\le g_\alpha\ast h\ \text{ on }[0,T].
\]
\end{lemma}
\begin{proof}
Since $g_\alpha\ast g_{1-\alpha}=1$,
we have $\Phi=1\ast \Phi'\le g_\alpha\ast h$.
\end{proof}

\begin{lemma}[$L^1$ theory for regular solutions]\label{lem:wp.basic}
Let $u$ and $v$ be regular solutions on $[0,T]$ with the same flux and data $u_0$ and $v_0$. Then for $0\le t\le T$ and $1\le p\le\infty$,
\[
\begin{gathered}
\int_\rr\big(u(t)-v(t)\big)_+\,dx\le\int_\rr(u_0-v_0)_+\,dx,\qquad \|u(t)-v(t)\|_{L^1(\rr)}\le\|u_0-v_0\|_{L^1(\rr)},\\
\|u(t)\|_{L^p(\rr)}\le\|u_0\|_{L^p(\rr)},\qquad \int_\rr u(t)\,dx=\int_\rr u_0\,dx,\\
\min\{0,\inf u_0\}\le u(t)\le\max\{0,\sup u_0\}.
\end{gathered}
\]
In particular, $u_0\ge0$ implies $u\ge0$.
\end{lemma}

\begin{proof}
Let $\phi=u-w$, with $w=v$ or $w\equiv \ell\ge0$. By (\ref{R3}) and Lemmas~\ref{lem:kato}, \ref{lem:dissip} and \ref{lem:L3},
\[
g_{1-\alpha}\ast\Phi_+'\le0\quad\text{on }(0,T],\qquad\text{hence}\qquad \Phi_+(t)\le\Phi_+(0).
\]
For $w=v$ this is the first bound, and exchanging $u$ and $v$ gives the contraction.
For $w\equiv \ell$ we apply it to $u$ and to $-u$, which is a regular solution with flux $-f(-\sigma)$. This gives
\[
\int_\rr(\pm u(t)-\ell)_+\,dx\le\int_\rr(\pm u_0-\ell)_+\,dx,\qquad \ell\ge0.
\]
With $\ell=\max\{0,\sup u_0\}$, respectively $\ell=\max\{0,-\inf u_0\}$, and the continuity of $u$, we obtain the pointwise bounds. Adding the two inequalities gives
\[
\int_\rr(|u(t)|-\ell)_+\,dx\le\int_\rr(|u_0|-\ell)_+\,dx,\qquad \ell\ge0.
\]
For $\ell=0$ this is the case $p=1$, and the cases $1<p<\infty$ follow from
\[
r^p=p(p-1)\int_0^\infty \ell^{p-2}(r-\ell)_+\,d\ell,\qquad r\ge0.
\]
Finally, by (\ref{R2})--(\ref{R3}), and since $u_x$ and $f(u)$ tend to $0$ at $\pm\infty$,
\[
g_{1-\alpha}\ast\frac{d}{dt}\int_\rr u\,dx=\int_\rr A(u)\,dx=0 .
\]
Lemma~\ref{lem:L3}, applied to $\pm\big(\int_\rr u(t)\,dx-\int_\rr u_0\,dx\big)$, gives the conservation of mass.
\end{proof}

The next remark explains how the estimates for regular solutions are transferred to mild solutions.

\begin{remark}[Passage to the limit]\label{rem:transfer}
Each estimate of Theorems~\ref{thm:L1} and \ref{thm:flux} is proved for regular solutions with an arbitrary regular flux. It depends on the data only through norms $\|u_0\|_{L^p(\rr)}$, in which it is nondecreasing.
It passes to the mild solution by the standard approximation argument: we apply it to the regular solutions $v_n$ of Theorem~\ref{thm:wp}(iii), with an admissible approximation of the flux and the standard data $u_{0,n}$, for which $\|u_{0,n}\|_{L^p(\rr)}\le\|u_0\|_{L^p(\rr)}$, and we let $n\to\infty$.
Bounds in $L^p(\rr)$ with $p<\infty$, such as those on the increments of the solution, pass to the limit by Theorem~\ref{thm:wp}(iii).
We use this in the proof of Theorem~\ref{thm:intro.L1}, for the estimate \eqref{eq:L1.incr} of Theorem~\ref{thm:L1}.
The bound \eqref{eq:L1.op} involves derivatives of $u$, and it passes to the limit in Remark~\ref{rem:regularity}.

Assume moreover that, for some $b\ge0$,
\[
|f_n(v_n(t,x))-\partial_xv_n(t,x)|\le b\qquad\text{for all }x\in\rr\text{ and }n\ge1.
\]
Then the functions $v_n(t)$ are bounded by $R$ and, by \eqref{eq:adm}, they are $(b+L_f(R)R)$-Lipschitz. By the Arzel\`a--Ascoli theorem every subsequence has a locally uniformly convergent subsequence,
and the $L^1$ convergence identifies the limit. Hence $v_n(t)$ converges locally uniformly to a Lipschitz representative $\tilde u$ of $u(t)$.
By \eqref{eq:adm}, $f_n(v_n(t))\to f(\tilde u)$ locally uniformly, and $\partial_xv_n(t)\to\tilde u'$ weakly-$*$ in $L^\infty(\rr)$.
Since balls of $L^\infty(\rr)$ are weakly-$*$ closed, $|f(\tilde u)-\tilde u'|\le b$ a.e.
We use this, with the flux bound of Theorem~\ref{thm:flux}, in the proofs of Corollary~\ref{cor:decay} and Theorem~\ref{thm:intro.L1}.
\end{remark}

The last lemma of this subsection is the equation tested with a fixed function of $x$. It is used for the lower bounds of Section~\ref{sec:L1}.

\begin{lemma}[Tested Volterra identity]\label{lem:tested}
Let $u$ be the mild solution with flux $f$ and datum $u_0\in L^1\cap L^\infty(\rr)$. For every $\zeta\in C_c^\infty(\rr)$ and $t\ge0$,
\[
\int_\rr\big(u(t)-u_0\big)\zeta\,dx=\int_0^tg_\alpha(t-s)\Big[\int_\rr u(s)\zeta''\,dx+\int_\rr f(u(s))\zeta'\,dx\Big]ds .
\]
\end{lemma}
\begin{proof}
We multiply the identity of Lemma~\ref{lem:A-fourier} by $\overline{\hat\zeta(\xi)}$ and integrate in $\xi$. Fubini's theorem and Plancherel's identity give the claim.
\end{proof}

\subsection{Weak solutions of \eqref{eq:general}}\label{sec:wp.general}
We consider \eqref{eq:general}, where $0<\lambda\le\Lambda$ and, for every $R>0$, $L_R$ is a constant, under the following assumptions.
\begin{itemize}
\item[(A)]\namedlabel{hyp:A}{A} $a$ is measurable on $(0,\infty)\times\rr$ and $\lambda\le a\le\Lambda$.
\item[(F)]\namedlabel{hyp:F}{F} $F$ is measurable in $t$. For a.e.\ $t$, $F(t,\cdot,\cdot)\in C^1(\rr^2)$, $F(t,x,0)=0$, and for $|\sigma|\le R$
\[
|\partial_\sigma F(t,x,\sigma)|\le L_R,\qquad |\partial_xF(t,x,\sigma)|\le L_R|\sigma|,\qquad \sigma\,\partial_xF(t,x,\sigma)\ge0 .
\]
\item[(Q)]\namedlabel{hyp:Q}{Q} $Q$ is measurable in $(t,x)$, $Q(t,x,0)=0$, and for $|\sigma|,|\sigma'|\le R$
\[
|Q(t,x,\sigma)-Q(t,x,\sigma')|\le L_R|\sigma-\sigma'|,\qquad \sigma\,Q(t,x,\sigma)\le0 .
\]
\end{itemize}
Condition (\ref{hyp:F}) holds for every flux $F(\sigma)$, in particular for the power flux, and (\ref{hyp:Q}) holds for $-c(t,x)|\sigma|^{q-1}\sigma$ with $c\ge0$ bounded.

\begin{definition}\label{def:weak}
Let $u_0\in L^1\cap L^\infty(\rr)$. A \emph{weak solution} of \eqref{eq:general} is a function $u\in L^\infty((0,\infty)\times\rr)$ with $u\in L^2(0,T;H^1(\rr))$ for every $T>0$, such that
\[
\int_0^T\!\!\int_\rr\Big(-\big(g_{1-\alpha}\ast(u-u_0)\big)\,\partial_t\varphi+a\,u_x\varphi_x-F(t,x,u)\,\varphi_x-Q(t,x,u)\,\varphi\Big)\,dx\,dt=0
\]
for every $T>0$ and every $\varphi\in C^1([0,T];H^1(\rr))$ with $\varphi(T)=0$.
\end{definition}

The identity says that $g_{1-\alpha}\ast(u-u_0)\in H^1(0,T;H^{-1}(\rr))$, with value $0$ at $t=0$ and time derivative equal to 
\[
\partial_x(au_x)-\big(F(t,x,u)\big)_x+Q(t,x,u).
\]
For $F=Q=0$, existence and uniqueness of weak solutions, also without the condition $u\in L^\infty$, are proved in \cite[Theorem 3.1]{Zacher09}.
For a.e.\ $t$ let $\mathcal A(t):H^1(\rr)\to H^{-1}(\rr)$ be given by
\[
\langle\mathcal A(t)v,\varphi\rangle=\int_\rr\big(a\,v_x\varphi_x-F(t,x,v)\varphi_x-Q(t,x,v)\varphi\big)\,dx ,
\]
so that \eqref{eq:general} reads
\[
\partial_t\big(g_{1-\alpha}\ast(u-u_0)\big)+\mathcal A(t)u=0 .
\]
For \eqref{convection} we have $\mathcal A(t)u=-A(u)$.
The sign conditions in (\ref{hyp:F}) and (\ref{hyp:Q}) make $\mathcal A(t)$ coercive on the test functions of the energy method.

\begin{lemma}\label{lem:coercive}
Assume \textup{(\ref{hyp:A})}, \textup{(\ref{hyp:F})} and \textup{(\ref{hyp:Q})}. For a.e.\ $t$, every $w\in H^1\cap L^\infty(\rr)$ and $p\ge2$,
\[
\langle\mathcal A(t)w,|w|^{p-2}w\rangle\ge\lambda(p-1)\int_\rr|w|^{p-2}w_x^2\,dx .
\]
\end{lemma}

\begin{proof}
Let $L=L_R$ with $R=\|w\|_{L^\infty(\rr)}$. The function
\[
\mathcal F_p(x,s)=\int_0^s|\sigma|^{p-2}F(t,x,\sigma)\,d\sigma
\]
satisfies $\partial_x\mathcal F_p\ge0$, by (\ref{hyp:F}), and $\mathcal F_p(\cdot,w)\in W^{1,1}(\rr)$, since $|F|+|\partial_xF|\le2L|\sigma|$ for $|\sigma|\le R$. The derivative of $x\mapsto\mathcal F_p(x,w(x))$ has integral zero, that is,
\[
\int_\rr(\partial_x\mathcal F_p)(x,w)\,dx+\int_\rr|w|^{p-2}F(t,x,w)\,w_x\,dx=0 .
\]
Hence
\[
-\int_\rr F(t,x,w)\big(|w|^{p-2}w\big)_x\,dx=(p-1)\int_\rr(\partial_x\mathcal F_p)(x,w)\,dx\ge0,
\]
and $-\int_\rr Q(t,x,w)|w|^{p-2}w\,dx\ge0$ by (\ref{hyp:Q}). Since $(|w|^{p-2}w)_x=(p-1)|w|^{p-2}w_x$ and $a\ge\lambda$, this gives the claim.
\end{proof}

\begin{theorem}\label{thm:wp.general}
Assume \textup{(\ref{hyp:A})}, \textup{(\ref{hyp:F})} and \textup{(\ref{hyp:Q})}, and let $u_0\in L^1\cap L^\infty(\rr)$. Then \eqref{eq:general} has a unique weak solution $u$, and
\[
|u|\le\|u_0\|_{L^\infty(\rr)},\qquad \|u(t)\|_{L^1(\rr)}\le M\quad\text{for a.e.\ }t>0.
\]
Moreover, let $T>0$ and let $k_n$ be as in Lemma~\ref{lem:yosida}. For $n$ large there is $u_n\in L^\infty(0,T;L^2(\rr))$, with $u_n(t)\in H^1(\rr)$ for a.e.\ $t$, such that
\[
|u_n|\le\|u_0\|_{L^\infty(\rr)},\qquad\|u_n(t)\|_{L^1(\rr)}\le M\quad\text{for a.e.\ }t,
\]
and
\begin{equation}\label{eq:reg}
\partial_t\big(k_n\ast(u_n-u_0)\big)+\mathcal A(t)u_n=0\quad\text{in }H^{-1}(\rr),\ \text{for a.e.\ }t\in(0,T),
\end{equation}
and a subsequence of $(u_n)$ converges to $u$ a.e.\ in $(0,T)\times\rr$. For \eqref{convection}, the mild solution of Theorem~\ref{thm:wp}(ii) is the weak solution.
\end{theorem}

\section{The linear equation}\label{sec:linear}

In this section we recall the classical energy method for the linear equation and we prove  the correct decay by a dyadic version of it.
Throughout the section $u$ solves
\[
\Da(u-u_0)-u_{xx}=0,\qquad u_0\in L^1\cap L^2(\rr).
\]
Here $u(t)=S(t)\star u_0$, with $S(t)$ the kernel of Section~\ref{sec:kernels}, so that $\|u(t)\|_{L^1(\rr)}\le M$ by (\ref{K1}).
When we apply \eqref{eq:I} to $u$, we first take $u_0\in C_c^\infty(\rr)$, for which $u$ is a regular solution with $f=0$ (Theorem~\ref{thm:wp}(i)).
We then pass to the limit in the final bounds, which involve only $\|u_0\|_{L^1(\rr)}$ and $\|u_0\|_{L^2(\rr)}$, using $\|S(t)\star g\|_{L^p(\rr)}\le\|g\|_{L^p(\rr)}$.
Multiplying \eqref{liniar} by $u$ and integrating gives
\begin{equation}\label{energy.lin}
\int_\rr u\,\Da(u-u_0)\,dx+\int_\rr u_x^2\,dx=0 .
\end{equation}

As explained in the introduction, this identity, with the Cauchy--Schwarz and Nash inequalities, gives the decay rate $ t^{-\alpha/5}$,
which is the rate of \cite[Theorem 6.1]{Kemppainen2014} for $d=1$.
We now adapt this method to obtain the optimal decay rate. 
For $j\in\zz$ let
\[
\mathcal C_j=\{\xi\in\rr:\ 2^{-j}\le|\xi|<2^{1-j}\},\qquad |\mathcal C_j|=2^{1-j},
\]
and define $\Delta_j$ by $\widehat{\Delta_j g}=1_{\mathcal C_j}\hat g$.
The operators $\Delta_j$ are orthogonal projections commuting with $\partial_x$ and $\Da$, and
\begin{align}
\label{eq:D1}&\|g\|^2_{L^2(\rr)}=\sum_{j\in\zz}\|\Delta_jg\|^2_{L^2(\rr)},\\
\label{eq:D2}&\|\Delta_jg\|^2_{L^2(\rr)}\le |\mathcal C_j|\,\|\hat g\|^2_{L^\infty(\rr)}\le \frac{2^{-j}}{\pi}\|g\|_{L^1(\rr)}^2,\\
\label{eq:D3}&2^{-j}\|\Delta_jg\|_{L^2(\rr)}\le\|\partial_x\Delta_jg\|_{L^2(\rr)}\le 2^{1-j}\|\Delta_jg\|_{L^2(\rr)}.
\end{align}

The point of the proof below is that it uses only the energy inequality \eqref{eq:I} on each block and the bound \eqref{eq:D2} on the data.
 
Let $y_j(t)=\|\Delta_ju(t)\|_{L^2(\rr)}$. Each block $\Delta_ju$ solves \eqref{liniar} with initial datum $\Delta_ju_0$, and it satisfies the energy identity
\[
\int _{\rr} \Delta_ju\, \Da(\Delta_ju-\Delta_ju_0)\,dx +\int _{\rr}(\Delta_ju)_{x}^{2}\,dx=0.
\]
By \eqref{eq:I} and \eqref{eq:D3},
\[
y_j\,\Da\big(y_j-y_j(0)\big)+2^{-2j}y_j^2\le0.
\]
If $y_j(0)=0$, then $y_j\equiv0$. Otherwise $y_j>0$, since the multiplier $E_\alpha(-\xi^2t^\alpha)$ is positive, and dividing by $y_j$,
\[
\Da\big(y_j-y_j(0)\big)+2^{-2j}y_j\le0.
\]
Lemma~\ref{lem:comparison} with $\nu=2^{-2j}$ and $h=0$ gives $y_j(t)\le y_j(0)E_\alpha(-2^{-2j}t^\alpha)$. Since $0<E_\alpha\le1$, \eqref{eq:ML} and \eqref{eq:D2} give
\[
y_j(t)^2\le y_j(0)^2E_\alpha(-2^{-2j}t^\alpha)\le \frac{M^2}{\pi}\,\frac{2^{-j}}{1+2^{-2j}t^\alpha}.
\]
For every $s>0$ we have
\[
\sum_{j\in\zz}\frac{2^{-j}}{1+2^{-2j}s}\le 5\,s^{-1/2}.
\]
Summing over $j\in\zz$, by \eqref{eq:D1} and this bound with $s=t^\alpha$, we obtain \eqref{eq:decay.liniar.1}.

\section{Energy methods}\label{sec:energy.nonlinear}

In this section we prove Theorem~\ref{thm:intro.energy}. The flux and reaction terms of \eqref{eq:general} have a sign (Lemma~\ref{lem:coercive}), and for \eqref{convection} the convection drops out.
Unlike the classical energy method presented in the introduction, we keep the whole identity \eqref{eq:Rk}. Since the weak solutions need not be regular in time, we use it with the Yosida kernels of Section~\ref{sec:yosida}.
The data then enter only through the term
\[
k_n(t)\int_\rr|u|^{p-2}u\,u_0\,dx,
\]
which is controlled by the mass. The time derivative appears through
\[
\cD_{k_n}\big[\|u\|_{L^p(\rr)}\big]=\partial_t\big(k_n\ast\|u\|_{L^p(\rr)}\big),
\]
an approximation of the Riemann--Liouville derivative of the norm, which does not contain the datum.

\subsection{The \texorpdfstring{$L^p$}{Lp} inequality}
The first step is the energy inequality for the Yosida approximations of Theorem~\ref{thm:wp.general}.

\begin{lemma}\label{lem:LPgeneral}
Assume \textup{(\ref{hyp:A})}, \textup{(\ref{hyp:F})} and \textup{(\ref{hyp:Q})}, let $p\ge2$, and let $N_p(t)=\|u_n(t)\|_{L^p(\rr)}$ with $u_n$ as in Theorem~\ref{thm:wp.general}. Then, for a.e.\ $t\in(0,T)$,
\[
N_p(t)^{p-1}\cD_{k_n}[N_p](t)+\lambda(p-1)\int_\rr|u_n(t)|^{p-2}\big(\partial_xu_n(t)\big)^2dx\le k_n(t)\int_\rr|u_n(t)|^{p-2}u_n(t)\,u_0\,dx .
\]
\end{lemma}

\begin{proof}
Since $k_n(t)+\int_0^t(-k_n'(t-s))\,ds=k_n(0)=n$, identity \eqref{eq:Rk} shows that \eqref{eq:reg} reads
\[
nu_n(t)+\mathcal A(t)u_n(t)=k_n(t)u_0+\int_0^tu_n(s)\big(-k_n'(t-s)\big)\,ds .
\]
We test it with $|u_n(t)|^{p-2}u_n(t)$, and we use Lemma~\ref{lem:coercive} and H\"older's inequality
\[
\int_\rr|u_n(t)|^{p-2}u_n(t)\,u_n(s)\,dx\le N_p(t)^{p-1}N_p(s).
\]
Since $-k_n'\ge0$, we obtain
\begin{align*}
&nN_p(t)^p+\lambda(p-1)\int_\rr|u_n(t)|^{p-2}\big(\partial_xu_n(t)\big)^2dx\\
&\qquad\le k_n(t)\int_\rr|u_n(t)|^{p-2}u_n(t)\,u_0\,dx+N_p(t)^{p-1}\int_0^tN_p(s)\big(-k_n'(t-s)\big)\,ds .
\end{align*}
Since, by \eqref{eq:Rk},
\[
nN_p(t)-\int_0^tN_p(s)\big(-k_n'(t-s)\big)\,ds=\cD_{k_n}[N_p](t),
\]
this gives the desired result.
\end{proof}

\subsection{Comparison with a barrier and decay}
The key estimate is a decay lemma for an abstract function that satisfies the inequalities of Lemma~\ref{lem:LPgeneral}, with the datum term bounded by the mass. We state it for a general kernel, so that it applies to all the Yosida kernels with the same constant.

\begin{lemma}\label{lem:abstract}
Let $T>0$, $\Pi\ge1$, $M>0$, and let $k\in W^{1,1}(0,T)$ be nonnegative and nonincreasing, with
\[
k\le\Pi g_{1-\alpha},\qquad -k'\le4\Pi|g_{1-\alpha}'|\qquad\text{on }(0,T).
\]
Let $u:(0,T)\to L^2(\rr)$ be strongly measurable and bounded in $L^1\cap L^\infty(\rr)$, and let $N_p(t)=\|u(t)\|_{L^p(\rr)}$. Assume that for a.e.\ $t$, $u(t)\in H^1(\rr)$, $N_1(t)\le M$ and, for $p=2^m$, $m\ge1$,
\begin{equation}\label{eq:abstractLP}
N_p(t)^{p-1}\cD_k[N_p](t)+(p-1)\int_\rr|u(t)|^{p-2}u_x^2(t)\,dx\le k(t)M\|u(t)\|_{L^\infty(\rr)}^{p-1}.
\end{equation}
Then
\[
\|u(t)\|_{L^\infty(\rr)}\le CMt^{-\alpha/2}\qquad\text{for a.e.\ }t\in(0,T),
\]
where $C$ depends only on $\alpha$ and $\Pi$.
\end{lemma}

\begin{proof}
The proof is an induction on $p=2^m$. For each $p$ we derive from \eqref{eq:abstractLP} an inequality that links two consecutive norms (Step 1), we bound from below the memory term of a barrier (Step 2), and we compare $N_p$ with a barrier (Step 3). Step 4 controls the constants.
Let $p=2^m$, $m\ge1$, and
\[
v=|u(t)|^{p/2},\qquad D=\|v_x\|_{L^2(\rr)},\qquad\text{so that}\qquad (p-1)\int_\rr|u|^{p-2}u_x^2\,dx=\frac{4(p-1)}{p^2}D^2 .
\]

\emph{Step 1. Absorption of the datum term.}
By   inequality \eqref{agmon} for $v$,
\[
\|u\|_{L^\infty(\rr)}^{p-1}\le N_p^{\frac{p-1}2}D^{\frac{p-1}p},
\]
and the right-hand side of \eqref{eq:abstractLP} is at most $bD^{\frac{p-1}p}$, with $b=\Pi g_{1-\alpha}(t)MN_p^{(p-1)/2}$.
We keep half of the term with $D^2$ on the left-hand side and use the other half to absorb the right-hand side. For $b>0$,
\[
\sup_{D>0}\Big(bD^{\frac{p-1}p}-\frac{2(p-1)}{p^2} D^2\Big)=\frac{2(p+1)}{p^2}\Big(\frac{bp}4\Big)^{\frac{2p}{p+1}},
\]
the supremum being attained at $D^{\frac{p+1}p}=\frac{bp}4$. This and Nash's inequality \eqref{nash} for $v$ give, after dividing by $N_p^{p-1}$, for a.e.\ $t$ with $N_p(t)>0$,
\begin{equation}\label{eq:star}
\cD_k[N_p]+\frac{2(p-1)}{p^2}\frac{N_p^{2p+1}}{N_{p/2}^{2p}}\le\frac{2(p+1)}{p^2}\Big(\frac p4\Pi g_{1-\alpha}(t)M\Big)^{\frac{2p}{p+1}}N_p^{-\frac{p-1}{p+1}} .
\end{equation}


\emph{Step 2. The memory term of a barrier.}
Let $0<\theta\le\frac12$ and $K>0$, and consider the barrier
\[
B(s)=Ks^{-\theta},\qquad s>0.
\]
We prove that it satisfies
\begin{equation}
  \label{eq:barier.1}
\cD_k[B](t)\ge-c_BKt^{-\theta-\alpha},\qquad c_B:=\frac{4\Pi\alpha\sqrt\pi}{\Gamma(\frac32-\alpha)} .
\end{equation}

By the bound on $k'$, the right-hand side of \eqref{eq:Rk} defines $\cD_k[B](t)$ for every $t\in(0,T)$.
We drop $k(t)B(t)\ge0$, and we use that $B$ is decreasing, the bound on $k'$ and the substitution $s=rt$:
\[
\cD_k[B](t)\ge-4\Pi\!\int_0^t\!\big(B(s)-B(t)\big)|g_{1-\alpha}'(t-s)|\,ds=-\frac{4\Pi\alpha K}{\Gamma(1-\alpha)}\,t^{-\theta-\alpha}\!\int_0^1\!\frac{r^{-\theta}-1}{(1-r)^{1+\alpha}}\,dr .
\]
Since $\theta\le\frac12$, for $0<r<1$ we have
\[
0\le r^{-\theta}-1\le r^{-1/2}-1\le r^{-1/2}(1-r),
\]
and the last integral is at most
\[
\int_0^1r^{-1/2}(1-r)^{-\alpha}\,dr=\frac{\sqrt\pi\,\Gamma(1-\alpha)}{\Gamma(\frac32-\alpha)} .
\]
This gives the desired result.

\emph{Step 3. Comparison with the barrier.}
Let $\theta_p=\frac\alpha2(1-\frac1p)$,
\[
c_p=\frac{2(p+1)}{p^2}\Big(\frac{p\Pi}4\Big)^{\frac{2p}{p+1}},
\]
and define $K_p$ by $K_1=M$ and
\[
K_p=\max\Big\{1,\Big(\frac{p^2(c_B+c_p)}{2(p-1)}\Big)^{\frac1{2p}}\Big\}K_{p/2},\qquad p=2^m,\ m\ge1.
\]
We prove by induction that
\[
N_p(t)\le K_pt^{-\theta_p}\qquad\text{for a.e.\ }t\in(0,T).
\]
We may assume $k(0)>0$, since otherwise $k\equiv0$, and \eqref{eq:abstractLP} with $p=2$ gives $u=0$. For $p=1$ the claim is $N_1\le M$. Let $p=2^m$, $m\ge1$, assume the claim for $p/2$, and let
\[
\theta=\theta_p, \qquad B(t)=K_pt^{-\theta}.
\]
Note that $K_p\ge K_{p/2}\ge M$.
  By \eqref{eq:star}, the induction hypothesis and  the bound $g_{1-\alpha}(t)=t^{-\alpha}/\Gamma(1-\alpha)\leq t^{-\alpha}$, for a.e.\ $t$ with $N_p(t)>0$ we have
\[
\cD_k[N_p](t)+G\big(t,N_p(t)\big)\le0,
\]
where, for $s>0$,
\[
G(t,s)=\frac{2(p-1)}{p^2}K_{p/2}^{-2p}\,t^{2p\theta_{p/2}}\,s^{2p+1}-c_p\big(Mt^{-\alpha}\big)^{\frac{2p}{p+1}}s^{-\frac{p-1}{p+1}} .
\]

We prove now that the barrier $B(t)=K_pt^{-\theta}$ is a supersolution. Since $2p\theta_p=\alpha(p-1)$ and $\theta_p-\theta_{p/2}=\frac\alpha{2p}$, all powers of $t$ in $G(t,B(t))$ are equal to $-\theta-\alpha$. Hence, by $M\le K_p$,
\[
G\big(t,B(t)\big)=K_pt^{-\theta-\alpha}\Big[\frac{2(p-1)}{p^2}\Big(\frac {K_p}{K_{p/2}}\Big)^{2p}-c_p\Big(\frac M{K_p}\Big)^{\frac{2p}{p+1}}\Big]\ge c_BK_pt^{-\theta-\alpha},
\]
and \eqref{eq:barier.1} gives
\[
\cD_k[B](t)+G\big(t,B(t)\big)\ge0,\qquad 0<t<T.
\]

We now compare. Since $k(0)=k(t)+\int_0^t(-k'(t-s))\,ds$, the right-hand side of \eqref{eq:Rk} gives,  {for a.e.\ $t$ with $N_p(t)>0$},
\begin{align*}
  G(t,B(t))-G(t,N_p(t))&\geq \cD_k[N_p](t)-\cD_k[B](t)\\
&=k(0)\big(N_p(t)-B(t)\big)-\int_0^t\big(N_p(s)-B(s)\big)\big(-k'(t-s)\big)\,ds.
\end{align*}

For a.e.\ $t$ with $N_p(t)>B(t)$ we have $N_p(t)>0$, and  {the left-hand side is nonpositive, since $G(t,\cdot)$ is increasing}. Hence, since $-k'\ge0$, the bounded function $w=(N_p-B)_+$ satisfies
\[
k(0)w(t)\le\int_0^tw(s)\big(-k'(t-s)\big)\,ds\qquad\text{for a.e.\ }t\in(0,T).
\]
Multiplying by $e^{-\mu t}$, $\mu>0$, we obtain
\[
k(0)\operatorname*{ess\,sup}_{t\in(0,T)}e^{-\mu t}w(t)\le\Big(\int_0^Te^{-\mu s}\big(-k'(s)\big)\,ds\Big)\operatorname*{ess\,sup}_{t\in(0,T)}e^{-\mu t}w(t).
\]
For $\mu$ large the integral is smaller than $k(0)$, so that $w=0$ a.e. Hence $N_p(t)\le K_pt^{-\theta_p}$ for a.e.\ $t$.

\emph{Step 4. The constants.}
The ratios $K_p/K_{p/2}\ge1$ do not depend on $M$, and $K_1=M$. Hence $K_{2^m}/M$ is nondecreasing in $m$, and $K_p\le CM$ for every $p$, with
\[
C=\lim_{m\to\infty}\frac{K_{2^m}}{M} .
\]
For $p\ge4$ we have $\frac{p\Pi}4\ge1$, so that
\[
c_p\le\frac{(p+1)\Pi^2}8,\qquad \frac{K_p}{K_{p/2}}\le\big((c_B+\Pi^2)p^2\big)^{1/(2p)} .
\]
Hence $C$ is finite and depends only on $c_B$ and $\Pi$, that is, on $\alpha$ and $\Pi$. Finally, for a.e.\ $t$,
\[
\|u(t)\|_{L^\infty(\rr)}=\lim_{p\to\infty}N_p(t)\le CMt^{-\alpha/2}. \qedhere
\]
\end{proof}

\subsection{Equations with rough coefficients}\label{sec:general}
We now apply Lemma~\ref{lem:abstract} to the Yosida approximations of the weak solutions of \eqref{eq:general}.

\begin{proof}[Proof of Theorem \ref{thm:intro.energy}]
We reduce to $\lambda=1$ by a dilation in $x$, and we apply Lemma~\ref{lem:abstract} to the Yosida approximations, uniformly in $n$.
Let $T>0$ and let $u_n$ be as in Theorem~\ref{thm:wp.general}. By Lemma~\ref{lem:LPgeneral} and the bound
\[
\int_\rr|u_n|^{p-2}u_n\,u_0\,dx\le M\|u_n\|_{L^\infty(\rr)}^{p-1},
\]
$u_n$ satisfies \eqref{eq:abstractLP} with $k=k_n$. By Lemma~\ref{lem:yosida}(a), Lemma~\ref{lem:abstract} applies with $\Pi=\Pi_\alpha$ and gives
\[
\|u_n(t)\|_{L^\infty(\rr)}\le CMt^{-\alpha/2}\qquad\text{for a.e.\ }t\in(0,T),
\]
with $C$ depending only on $\alpha$, hence independent of $n$ and $T$.
Along a subsequence, $u_n\to u$ a.e.\ in $(0,T)\times\rr$, hence $u_n(t)\to u(t)$ a.e.\ for a.e.\ $t$, and the bound passes to $u$.
\end{proof}

\begin{remark}[The convection--diffusion equation]\label{rem:parabolic}
The mild solution of \eqref{convection} is the weak solution of \eqref{eq:general} with $a\equiv1$, $Q=0$ and $F=f$ (Theorem~\ref{thm:wp.general}), so Theorem~\ref{thm:intro.energy} applies to it, with $\lambda=1$.
Since the mild solution is continuous with values in $L^1(\rr)$, the bound holds for every $t>0$: if it holds along $t_j\to t$, a subsequence of $u(t_j)$ converges to $u(t)$ a.e.
For $q\ge2$ the rate is optimal (Theorem~\ref{thm:lower}), and for $1<q<2$ the convection gives a faster decay (Section~\ref{sec:L1}).
\end{remark}

\section{\texorpdfstring{$L^1$}{L1} time regularity, a flux bound and sharp decay}\label{sec:L1}

The energy methods of Section~\ref{sec:energy.nonlinear} give the parabolic rate $t^{-\alpha/2}$ for every $q>1$. They do not see the convection, and by Theorem~\ref{thm:lower} below this rate is optimal for $q\ge2$.
For $1<q<2$ the convection gives a faster decay, and for large times the estimates of this section improve those of Section~\ref{sec:energy.nonlinear}.
They use the identity \eqref{acv2} in $L^1(\rr)$, together with the time scaling of the Caputo derivative and integration in $x$. They give the hyperbolic rate $t^{-\alpha/q}$ (Corollary~\ref{cor:decay}).

Up to Theorem~\ref{thm:flux}, $u$ is a regular solution on $[0,T]$ (Definition~\ref{def:regular}) with a regular flux $f$, and $A(u)=u_{xx}-f(u)_x$ as in Section~\ref{sec:L1theory}. By (\ref{R3}),
\[
(g_{1-\alpha}\ast\partial_tu)(t)=A(u(t)),\qquad t\in(0,T].
\]
No other property of $f$ is used, and the estimates depend only on $\alpha$ and $M$. They therefore pass to the mild solutions of \eqref{convection}.

We first explain the idea. Integrating $A(u)$ over $(-\infty,x)$ and using $\int_\rr A(u)\,dx=0$ gives
\[
|f(u)-u_x|(t,x)=\Big|\int_{-\infty}^xA(u(t))\,dy\Big|\le\tfrac12\|A(u(t))\|_{L^1(\rr)}.
\]
By the identity \eqref{acv2},
\[
\Gamma(1-\alpha)A(u(t))=\frac{u(t)-u_0}{t^\alpha}+\alpha\int_0^t\frac{u(t)-u(s)}{(t-s)^{1+\alpha}}\,ds .
\]
In $L^1(\rr)$ the datum term and the increments are bounded by $2M$, with no factor $\|u_0\|_{L^\infty(\rr)}$: in $L^1$ the memory of the data is the memory of the mass.
Since the kernel is not integrable at $s=t$, we also need the increments to be small near the diagonal. Any power modulus $M(1-s/t)^\beta$ with $\beta>\alpha$ would do, and the estimate \eqref{eq:L1.incr} of Theorem~\ref{thm:L1} gives the scale-invariant bound
\[
\|u(t)-u(s)\|_{L^1(\rr)}\lesssim M\Big(1-\big(\tfrac st\big)^\alpha\Big)\le M\Big(1-\frac st\Big),
\]
which is the case $\beta=1$.

Time translation, the usual route to time regularity, is not available: the memory of $\Da$ starts at $t=0$, so a translate of $u$ in time does not solve \eqref{convection}.
Dilations in time respect this memory: by Lemma~\ref{lem:scaling}, $v=u(\kappa\,\cdot)$ has the same datum as $u$ and solves
\[
g_{1-\alpha}\ast\partial_tv=\kappa^\alpha A(v).
\]
We compare $v$ with $u$. The $L^1$ dissipativity of $A$ removes $A(v)-A(u)$, and only the term $(\kappa^\alpha-1)A(v)$ remains. Hence the increments of $u$ are controlled by $A(u)$: if $\|A(u(t))\|_{L^1(\rr)}\le\Xi\,g_{1-\alpha}(t)$, the homogeneity $g_{1-\alpha}(\kappa t)=\kappa^{-\alpha}g_{1-\alpha}(t)$ gives $\|u(\tau)-u(s)\|_{L^1(\rr)}\le\Xi\,(1-(s/\tau)^\alpha)$.
Conversely, by \eqref{acv2}, $\Xi$ is controlled by the increments. Since they are bounded by $2M$ far from the diagonal, they enter only through their power $\alpha<1$. This closes the estimate.

\subsection{Time scaling}
The proof of Theorem~\ref{thm:L1} combines the Kato inequality (Lemma~\ref{lem:kato}), the dissipativity of $A$ (Lemma~\ref{lem:dissip}) and Lemma~\ref{lem:L3} with the following symmetry of $g_{1-\alpha}\ast\partial_t$.

\begin{lemma}[Time scaling]\label{lem:scaling}
Let $u$ be a regular solution on $[0,T]$ and $\kappa\ge1$. Then $v(t,x):=u(\kappa t,x)$ satisfies, for every $t\in(0,T/\kappa]$,
\[
(g_{1-\alpha}\ast\partial_tv)(t)=\kappa^\alpha A(v(t))\qquad\text{in }L^1(\rr).
\]
\end{lemma}

\subsection{Time regularity}
The following theorem is the key estimate of this section.

\begin{theorem}[$L^1$ time regularity]\label{thm:L1}
Let $\omega_\alpha$ be the unique root in $(1,\infty)$ of $\omega=1+\frac{\omega^\alpha}{1-\alpha}$. Every regular solution satisfies, for all $0\le s<\tau\le T$,
\begin{gather}
\|u(\tau)-u(s)\|_{L^1(\rr)}\le\min\Big\{2,\ 2\omega_\alpha\Big(1-\big(\tfrac s\tau\big)^\alpha\Big)\Big\}\,M,\label{eq:L1.incr}\\
\|A(u(\tau))\|_{L^1(\rr)}\le2\omega_\alpha M\,g_{1-\alpha}(\tau).\label{eq:L1.op}
\end{gather}
\end{theorem}

\begin{proof}
The proof consists of three steps. Step 1 bounds the operator by the increments, Step 2 compares $u$ with its time dilations, and Step 3 combines the two bounds into an inequality for the best constant $\Omega$ below.
  Assume $M>0$ and let
\[
\Omega:=\sup\Big\{\frac{\|u(\tau)-u(s)\|_{L^1(\rr)}}{1-(s/\tau)^\alpha}:\ 0\le s<\tau\le T\Big\}.
\]
By (\ref{R2}),
\[
\|u(\tau)-u(s)\|_{L^1(\rr)}\le \frac C\alpha(\tau^\alpha-s^\alpha)\le \frac{CT^\alpha}\alpha\Big(1-\big(\tfrac s\tau\big)^\alpha\Big),
\]
so that $\Omega<\infty$. By Lemma~\ref{lem:wp.basic}, $\|u(t)\|_{L^1(\rr)}\le M$, and therefore $\|u(\tau)-u(s)\|_{L^1(\rr)}\le2M$.

\emph{Step 1. The memory controls the operator.}
We apply \eqref{acv2} at every $x$, take the $L^1$ norm, and use Minkowski's inequality, the substitution $s=r\tau$ and the inequality $1-r^\alpha\le1-r$ on $[0,1]$:
\begin{align*}
\Gamma(1-\alpha)\|A(u(\tau))\|_{L^1(\rr)}
&\le\frac{\|u(\tau)-u_0\|_{L^1(\rr)}}{\tau^\alpha}+\alpha\int_0^\tau\frac{\|u(\tau)-u(s)\|_{L^1(\rr)}}{(\tau-s)^{1+\alpha}}\,ds\\
&\le\tau^{-\alpha}\Big[2M+\alpha\int_0^1\min\{2M,\Omega(1-r)\}(1-r)^{-1-\alpha}\,dr\Big].
\end{align*}
For $c>0$ and $c'\ge0$ we have
\[
\int_0^\infty\min\{c,c'r\}\,r^{-1-\alpha}\,dr=\frac{c^{1-\alpha}c'^{\,\alpha}}{\alpha(1-\alpha)} .
\]
Hence, with
\[
\Xi:=2M+\frac{(2M)^{1-\alpha}\Omega^\alpha}{1-\alpha},
\]
we obtain
\[
\|A(u(\tau))\|_{L^1(\rr)}\le \Xi\, g_{1-\alpha}(\tau),\qquad 0<\tau\le T.
\]
Using the mass bound far from the diagonal and $\Omega$ only near it makes $\Xi$ grow like $\Omega^\alpha$. This is the key point that allows us to close Step 3. 

\emph{Step 2. Comparison with the time dilations.}
Fix $\kappa>1$, let $v=u(\kappa\,\cdot)$ and $\phi=v-u$ on $[0,T/\kappa]$, so that $\phi(0)=0$. By Lemma~\ref{lem:scaling},
\[
g_{1-\alpha}\ast\partial_t\phi=\big[A(v)-A(u)\big]+(\kappa^\alpha-1)A(v).
\]
We multiply by $\sgn\phi$ and integrate in $x$. By Lemma~\ref{lem:dissip}, for $0<t\le T/\kappa$,
\begin{align*}
\int_\rr\sgn\phi\,(g_{1-\alpha}\ast\partial_t\phi)(t)\,dx
&=\int_\rr\sgn\phi\,\big[A(v(t))-A(u(t))\big]\,dx+(\kappa^\alpha-1)\int_\rr\sgn\phi\,A(v(t))\,dx\\
&\le(\kappa^\alpha-1)\|A(u(\kappa t))\|_{L^1(\rr)}.
\end{align*}
By Lemma~\ref{lem:kato} and Step 1, $\Phi=\|\phi\|_{L^1(\rr)}$ satisfies
\[
(g_{1-\alpha}\ast\Phi')(t)\le(\kappa^\alpha-1)\|A(u(\kappa t))\|_{L^1(\rr)}\le(\kappa^\alpha-1)\,\Xi\, g_{1-\alpha}(\kappa t)=(1-\kappa^{-\alpha})\,\Xi\, g_{1-\alpha}(t).
\]
Lemma~\ref{lem:L3} and $g_\alpha\ast g_{1-\alpha}=1$ give
\[
\|u(\kappa t)-u(t)\|_{L^1(\rr)}=\Phi(t)\le(1-\kappa^{-\alpha})\,\Xi,\qquad 0\le t\le T/\kappa.
\]

\emph{Step 3. Closing the estimate.}
For $0<s<\tau\le T$ we take $t=s$ and $\kappa=\tau/s$ in Step 2:
\[
\|u(\tau)-u(s)\|_{L^1(\rr)}\le\Big(1-\big(\tfrac s\tau\big)^\alpha\Big)\,\Xi.
\]
For $s=0$ this also holds, since $\|u(\tau)-u_0\|_{L^1(\rr)}\le2M\le \Xi$. Hence
\[
\Omega\le \Xi=2M+\frac{(2M)^{1-\alpha}\Omega^\alpha}{1-\alpha},\qquad\text{that is,}\qquad w\le1+\frac{w^\alpha}{1-\alpha},\quad w=\frac{\Omega}{2M},
\]
which proves $w\leq \omega_\alpha$.
Since $\Xi/(2M)=1+w^\alpha/(1-\alpha)$ is increasing in $w$,
\[
\Xi\le2M\Big(1+\frac{\omega_\alpha^\alpha}{1-\alpha}\Big)=2\omega_\alpha M,
\]
so that Step 3 gives \eqref{eq:L1.incr} and Step 1 gives \eqref{eq:L1.op}.
\end{proof}

\subsection{The flux bound}
Integrating the equation in $x$ turns \eqref{eq:L1.op} into a pointwise bound.

\begin{theorem}[Pointwise flux bound]\label{thm:flux}
Every regular solution satisfies, for all $t\in(0,T]$ and $x\in\rr$,
\[
\big|f(u(t,x))-u_x(t,x)\big|\le \omega_\alpha M\,g_{1-\alpha}(t)=\frac{\omega_\alpha M}{\Gamma(1-\alpha)\,t^\alpha}.
\]
\end{theorem}

\begin{proof}
Fix $t$. Since $u(t)\in W^{2,1}(\rr)$, $u_x$ and $f(u)$ vanish at infinity, so that $\int_\rr A(u(t))\,dy=0$ and
\[
f(u(t,x))-u_x(t,x)=-\int_{-\infty}^xA(u(t))\,dy=\int_x^\infty A(u(t))\,dy .
\]
Therefore $|f(u)-u_x|\le\frac12\|A(u(t))\|_{L^1(\rr)}$, and \eqref{eq:L1.op} gives the claim.
\end{proof}

\begin{remark}[Regularity of the mild solutions]\label{rem:regularity}
Let $u$ be the mild solution with flux $f$ and datum $u_0\in L^1\cap L^\infty(\rr)$, and let $v_n$ be the approximations of Remark~\ref{rem:transfer}, with fluxes $f_n$. For $t>0$, by \eqref{acv2},
\[
\Gamma(1-\alpha)\big(\partial_{xx}v_n-\partial_xf_n(v_n)\big)(t)=\frac{v_n(t)-u_{0,n}}{t^\alpha}+\alpha\int_0^t\frac{v_n(t)-v_n(s)}{(t-s)^{1+\alpha}}\,ds .
\]
By \eqref{eq:L1.incr}, the integrand has $L^1(\rr)$ norm at most $\min\{2M,2\omega_\alpha M(1-\frac st)\}(t-s)^{-1-\alpha}$, and dominated convergence shows that the right-hand side converges in $L^1(\rr)$. Since $v_n(t)\to u(t)$ and, by \eqref{eq:adm}, $f_n(v_n(t))\to f(u(t))$ in $L^1(\rr)$, the limit is the distribution $u_{xx}(t)-f(u(t))_x$, which we denote by $A(u(t))$. Hence $A(u(t))\in L^1(\rr)$, with $\|A(u(t))\|_{L^1(\rr)}\le2\omega_\alpha Mg_{1-\alpha}(t)$ by \eqref{eq:L1.op}, and $\partial_xv_n(t)-f_n(v_n(t))$ converges uniformly.  By Theorem~\ref{thm:flux} and Remark~\ref{rem:transfer}, $v_n(t)\to\tilde u$ and $f_n(v_n(t))\to f(\tilde u)$ locally uniformly, where $\tilde u$ is the Lipschitz representative of $u(t)$. Hence $\partial_xv_n(t)$ converges locally uniformly to a continuous function, and therefore $\tilde u\in C^1(\rr)$.
\end{remark}

\subsection{Decay}
Formally, at a maximum point of $|u(t)|$ we have $u_x=0$, and the flux bound gives $|u|^q\le\omega_\alpha Mg_{1-\alpha}(t)$. The diffusion term does not enter.
The power flux is not regular and the data are only in $L^1\cap L^\infty(\rr)$. We therefore argue on approximations, whose limit we treat as a Lipschitz function, and the following lemma replaces the condition $u_x=0$.

\begin{lemma}[Consequences of a flux bound]\label{lem:fluxbound}
Let $w:\rr\to\rr$ be Lipschitz with $w\in L^1(\rr)$.
\begin{enumerate}
\item[(a)] $w(x)\to0$ as $|x|\to\infty$, and
\[
\|w\|_{L^\infty(\rr)}^2\le\|w'\|_{L^\infty(\rr)}\|w\|_{L^1(\rr)}.
\]
\item[(b)] Let $f\in C(\rr)$ and $b\ge0$ be such that
\[
|f(w(x))-w'(x)|\le b\qquad\text{for a.e.\ }x\in\rr .
\]
If $w\not\equiv0$, then $|w|$ attains its maximum, and $|f(w(x_0))|\le b$ at every maximum point $x_0$ of $|w|$.
\item[(c)] If, in (b), $f(\sigma)=|\sigma|^{q-1}\sigma$ with $q>1$, then
\[
\|w\|_{L^\infty(\rr)}^q\le b,\qquad \|w'\|_{L^\infty(\rr)}\le2b,\qquad \|w\|_{L^\infty(\rr)}^2\le2b\|w\|_{L^1(\rr)}.
\]
\end{enumerate}
\end{lemma}

\begin{proof}
(b)  By (a), $|w|$ attains its maximum $m>0$ at some $x_0$, and   assume $w(x_0)=m$.
If $f(m)>b$, then by continuity $f(w)-b\ge\eta>0$ on some interval $[x_0,x_0+\delta]$, and
\[
w(x_0+\delta)\ge m+\int_{x_0}^{x_0+\delta}\big(f(w)-b\big)\,dx\ge m+\eta\delta ,
\]
a contradiction. If $f(m)<-b$, the same argument on an interval to the left of $x_0$ gives a contradiction.

(c) By (b), $\|w\|_{L^\infty(\rr)}^q=|f(w(x_0))|\le b$. Then $|w'|\le|f(w)|+b\le2b$ a.e., and (a) gives the last bound.
\end{proof}

\begin{corollary}[Decay for all data]\label{cor:decay}
Let $u$ be the mild solution of \eqref{convection} with datum $u_0\in L^1\cap L^\infty(\rr)$. Then, for every $t>0$, $u(t)$ has a Lipschitz representative, and:
\begin{enumerate}
\item[(i)] (hyperbolic bound)
\[
\|u(t)\|_{L^\infty(\rr)}^q\le \omega_\alpha Mg_{1-\alpha}(t);
\]
\item[(ii)] (parabolic bound)
\[
\|u(t)\|_{L^\infty(\rr)}\le M\big(2\omega_\alpha g_{1-\alpha}(t)\big)^{1/2}.
\]
\end{enumerate}
\end{corollary}

\begin{proof}
Let $v_n$ be the approximations of Remark~\ref{rem:transfer}. By Theorem~\ref{thm:flux},
\[
|f_n(v_n(t,x))-\partial_xv_n(t,x)|\le\omega_\alpha\|u_{0,n}\|_{L^1(\rr)}g_{1-\alpha}(t)\le\omega_\alpha Mg_{1-\alpha}(t).
\]
By Remark~\ref{rem:transfer}, $u(t)$ has a Lipschitz representative $\tilde u$ with
\[
|f(\tilde u)-\tilde u'|\le\omega_\alpha Mg_{1-\alpha}(t)\quad\text{a.e.},\qquad \|\tilde u\|_{L^1(\rr)}=\|u(t)\|_{L^1(\rr)}\le M .
\]
Lemma~\ref{lem:fluxbound}(c) with $w=\tilde u$ gives (i) and (ii).
\end{proof}

For large times, (i) for $q\le2$ and (ii) for $q\ge2$ are sharp for nonnegative data, as the following theorem shows.

\begin{theorem}[Lower bounds]\label{thm:lower}
Let $u$ be the mild solution of \eqref{convection} with datum $0\le u_0\in L^1\cap L^\infty(\rr)$, $u_0\neq0$, and let $p_*=\min\{q,2\}$. There are $\rho_0,c>0$, depending on $\alpha,q,u_0$, such that for all $t\ge1$ and $1\le p\le\infty$,
\[
\|u(t)\|_{L^p(\rr)}\ge c\,t^{-\frac{\alpha}{p_*}(1-\frac1p)},\qquad\text{in particular}\quad \|u(t)\|_{L^\infty(\rr)}\ge\frac{3M}{16\rho_0}\,t^{-\alpha/p_*}.
\]
\end{theorem}

\begin{proof}
We prove that, on the scale $\rho\simeq t^{\alpha/p_*}$, at most a quarter of the mass lies outside $(-2\rho,2\rho)$. This argument is inspired by the classical tail estimates for convection--diffusion problems  
\cite[Lemma 5]{laurencot2005asymptotic}, \cite{endal2023largetime}.

By Theorem~\ref{thm:wp}(ii), $u\ge0$ and $\int_\rr u(t)\,dx=\int_\rr u_0\,dx=M$. We fix a cut-off $\chi_0\in C^\infty(\rr)$ with
\[
0\le\chi_0\le1,\qquad \chi_0=0\ \text{ on }[-1,1],\qquad \chi_0=1\ \text{ outside }(-2,2),
\]
and in this proof $C$ denotes a constant depending only on $\alpha,q,M,\chi_0$.
Fix $t\ge1$ and $\rho>0$, and let
\[
\chi=\chi_0(\cdot/\rho),\qquad m(s)=\int_\rr\chi\, u(s)\,dx .
\]
By Lemma~\ref{lem:tested} with $\zeta=1-\chi\in C_c^\infty(\rr)$, and since $\int_\rr(u(s)-u_0)\,dx=0$,
\[
m(s)-m(0)=\Big(g_\alpha\ast\Big[\int_\rr\chi''u\,dx+\int_\rr\chi'f(u)\,dx\Big]\Big)(s)\le\Big(g_\alpha\ast\Big[\frac{CM}{\rho^2}+\frac{CM}{\rho}\|u(\cdot)\|_{L^\infty(\rr)}^{q-1}\Big]\Big)(s).
\]
By Corollary~\ref{cor:decay}(i), and since $\alpha-\alpha/q<1$,
\[
\|u(s)\|_{L^\infty(\rr)}^{q-1}\le Cs^{-\alpha+\alpha/q},\qquad \Big(g_\alpha\ast s^{-\alpha+\alpha/q}\Big)(t)=\frac{\Gamma(1-\alpha+\alpha/q)}{\Gamma(1+\alpha/q)}\,t^{\alpha/q}.
\]
Hence
\[
m(t)\le\int_{|x|>\rho}u_0\,dx+C\Big(\frac{t^\alpha}{\rho^2}+\frac{t^{\alpha/q}}{\rho}\Big).
\]
The two error terms come from diffusion, on the scale $t^{\alpha/2}$, and from convection, on the hyperbolic scale $t^{\alpha/q}$, controlled through the upper bound (i). For $t\ge1$ the larger scale is $t^{\alpha/p_*}$, and we take $\rho$ proportional to it:
\[
\rho=\rho_0t^{\alpha/p_*} .
\]
For $t\ge1$, since $p_*\le2$ and $p_*\le q$, this gives for {large} enough $\rho_0$  that
\[
m(t)\le\int_{|x|>\rho_0}u_0\,dx+C\big(\rho_0^{-2}+\rho_0^{-1} \big)\leq \frac M4.
\]
For such  $\rho_0$, by the conservation of mass, at least three quarters of the mass lies in $[-2\rho,2\rho]$, and H\"older's inequality on this interval gives
\[
\|u(t)\|_{L^\infty(\rr)}\ge\frac{3M}{16\rho},\qquad \|u(t)\|_{L^p(\rr)}\ge\frac34M(4\rho)^{-(1-\frac1p)},
\]
which gives the claim.
\end{proof}

\begin{proof}[Proof of Theorems~\ref{thm:intro.L1} and \ref{thm:intro.decay}]
The estimate \eqref{eq:L1.incr} and Theorem~\ref{thm:flux} hold for the approximations $v_n$ of Remark~\ref{rem:transfer}, whose data have $L^1$ norm at most $M$. They pass to $u$ by that remark, as in the proof of Corollary~\ref{cor:decay}.
The decay estimates are Corollary~\ref{cor:decay} and Theorem~\ref{thm:lower}.
For nonnegative data, the upper bound
\[
\|u(t)\|_{L^p(\rr)}\lesssim t^{-\frac{\alpha}{p_*}(1-\frac1p)}
\]
follows from (i) when $q<2$ and from (ii) when $q\ge2$, by interpolation with $\|u(t)\|_{L^1(\rr)}\le M$.
\end{proof}

\appendix

\section{Proofs of the results of Section~\ref{sec:wp.mild}}\label{app:wp}

In this appendix we sketch the proof of Theorem~\ref{thm:wp} and prove Lemma~\ref{lem:A-fourier}, which is used in Lemma~\ref{lem:tested}.
The arguments are standard for semilinear subdiffusion equations. We refer to \cite[Ch.~7, Section~7.1, p.~277]{Jin} for the fundamental solutions, and to \cite{LLW18,CNP15} for mild solutions with a nonlinearity in divergence form in $\rr^d$. We only indicate the main steps.

Let $\BUC(\rr)$ be the space of bounded uniformly continuous functions. We work in
\[
\EE=L^1\cap\BUC(\rr),\qquad\text{with}\qquad \|g\|_{L^1\cap\BUC}=\|g\|_{L^1(\rr)}+\|g\|_{L^\infty(\rr)} .
\]
On $L^1\cap\BUC(\rr)$, unlike on $L^1\cap L^\infty(\rr)$, translations and point evaluations are continuous. Moreover, $S(t)\star g\to g$ as $t\to0$, since $S(t)$ is an approximate identity by \eqref{eq:scaling}.
For $k\in\mathbb N$ let $\EE^k$ be the space of $C^k$ functions whose derivatives of order at most $k$ belong to $L^1\cap\BUC(\rr)$, with the norm
\[
\|g\|_{\EE^k}=\sum_{j\le k}\|g^{(j)}\|_{L^1\cap\BUC},\qquad\text{so that}\qquad W^{k+1,1}(\rr)\subset\EE^k\subset W^{k,1}(\rr).
\]
By (\ref{K1})--(\ref{K2}), for $g\in\EE$,
\[
\|S(t)\star g\|_\EE\le\|g\|_\EE,\qquad \|P(t)\star g\|_\EE\le g_\alpha(t)\|g\|_\EE,\qquad \|\partial_xP(t)\star g\|_\EE\le g_{\alpha/2}(t)\|g\|_\EE .
\]

We use the singular Gronwall inequality. Let $\beta>0$, $0\le\theta<1$ and $c,L\ge0$, and let $\varphi\ge0$ be measurable on $(0,T]$, with $\sup_tt^\theta\varphi(t)<\infty$ and
\[
\varphi(t)\le c\,t^{-\theta}+L\,(g_\beta\ast\varphi)(t),\qquad 0<t\le T.
\]
Then \cite[Ch.~4, Theorem~4.2(i)--(ii), pp.~98--99]{Jin}
\begin{equation}\label{eq:A-gronwall}
\varphi(t)\le c\,\Gamma(1-\theta)\,t^{-\theta}E_{\beta,1-\theta}(Lt^\beta),\qquad 0<t\le T.
\end{equation}

Accordingly, let $m\in C([0,T];\BUC(\rr))$ with $|m|\le L$, let $0\le\theta<1$, and let $\phi_0:(0,T]\to\EE$ be strongly measurable, with $t^\theta\|\phi_0(t)\|_\EE$ bounded. The linear Volterra equation
\begin{equation}\label{eq:A-lin}
\phi(t)=\phi_0(t)-\int_0^t\partial_xP(t-s)\star [m(s)\phi(s)]\,ds
\end{equation}
has a unique solution with $t^\theta\|\phi(t)\|_\EE$ bounded. It is given by its Neumann series and bounded by \eqref{eq:A-gronwall} with
\[
\beta=\frac\alpha2,\qquad c=\sup_{0<t\le T}t^\theta\|\phi_0(t)\|_\EE .
\]
It is continuous on $(0,T]$ if $\phi_0$ is. It is continuous on $[0,T]$ if moreover $\theta=0$ and $\phi_0\in C([0,T];\EE)$, since then the integral in \eqref{eq:A-lin} is $O(g_{1+\alpha/2}(t))$.

\subsection{Mild solutions}
The first lemma gives existence for Lipschitz fluxes and the stability estimate that drives all the approximations.

\begin{lemma}\label{lem:A-mild}
\textup{(a)} Let $f$ be Lipschitz with $f(0)=0$, and let $u_0\in\EE$. There is a unique $u\in C([0,\infty);\EE)$ that satisfies \eqref{eq:mild} for all $t\ge0$, with an $\EE$-valued Bochner integral.

\textup{(b)} For $i=1,2$ let $u_i$ be a mild solution on $[0,T]$ with flux $f_i$ and datum $u_0^i\in L^1\cap L^\infty(\rr)$. Let $R\ge\sup_{t\le T}\|u_i(t)\|_{L^\infty(\rr)}$, $i=1,2$, and
\[
L=\max_{|\sigma|\le R}|f_1'(\sigma)|,\qquad \delta=\max_{|\sigma|\le R}|f_1'(\sigma)-f_2'(\sigma)|,\qquad m_2=\max_{t\le T}\|u_2(t)\|_{L^1(\rr)}.
\]
Then for $0\le t\le T$,
\[
\|u_1(t)-u_2(t)\|_{L^1(\rr)}\le\big(\|u^1_0-u^2_0\|_{L^1(\rr)}+\delta\,m_2\,g_{1+\alpha/2}(T)\big)\,E_{\alpha/2}\big(Lt^{\alpha/2}\big).
\]
\end{lemma}

\begin{proof}
(a) By (\ref{K2}), the approximate identity property above and dominated convergence, the right-hand side of \eqref{eq:mild}, with an $\EE$-valued Bochner integral, maps $C([0,T];\EE)$ into itself.
Since $\|\partial_xP(t)\star g\|_\EE\le g_{\alpha/2}(t)\|g\|_\EE$, it is a contraction for the norm
\[
\sup_{t\le T}e^{-bt}\|u(t)\|_\EE
\]
with $b$ large. This is the argument of \cite[Ch.~6, proof of Theorem~6.17, p.~222]{Jin} (see \cite[Theorem 3.1]{JLZ18}). Uniqueness on every $[0,T]$ gives the global solution.

(b) Since $|u_i|\le R$ and $f_1(0)=f_2(0)=0$,
\[
|f_1(u_1)-f_2(u_2)|\le L|u_1-u_2|+\delta|u_2| .
\]
By (\ref{K1})--(\ref{K2}), the function $\varphi(t)=\|u_1(t)-u_2(t)\|_{L^1(\rr)}$ satisfies
\[
\varphi\le\|u_0^1-u_0^2\|_{L^1(\rr)}+\delta m_2g_{1+\alpha/2}(T)+L\,g_{\alpha/2}\ast\varphi,
\]
and \eqref{eq:A-gronwall} with $\theta=0$ gives the claim.
\end{proof}

The Fourier transform of \eqref{eq:mild} gives the Volterra form of the equation.

\begin{lemma}\label{lem:A-fourier}
Let $u$ be a mild solution on $[0,T]$ with flux $f$. Then for $0\le t\le T$ and $\xi\in\rr$
\[
\hat u(t,\xi)=\hat u_0(\xi)+\Big(g_\alpha\ast\big(-\xi^2\hat u(\cdot,\xi)-i\xi\,\widehat{f(u)}(\cdot,\xi)\big)\Big)(t),
\]
that is, $u(t)=u_0+\big(g_\alpha\ast(u_{xx}-f(u)_x)\big)(t)$ in the sense of distributions.
\end{lemma}

\begin{proof}
For fixed $\xi$ let $e_\xi(t)=E_\alpha(-\xi^2t^\alpha)$, and let $r_{\xi^2}(t)=t^{\alpha-1}E_{\alpha,\alpha}(-\xi^2t^\alpha)$ be the resolvent kernel of Section~\ref{sec:prelim}. Their power series and the identity $g_{\beta_1}\ast g_{\beta_2}=g_{\beta_1+\beta_2}$ give
\[
e_\xi=1-\xi^2g_\alpha\ast e_\xi,\qquad r_{\xi^2}=g_\alpha-\xi^2g_\alpha\ast r_{\xi^2} .
\]
The Fourier transform of \eqref{eq:mild} is
\[
\hat u(t)=e_\xi(t)\hat u_0-i\xi\,\big(r_{\xi^2}\ast\widehat{f(u)}\big)(t),
\]
and the two identities give the claim.
\end{proof}

\subsection{Regular solutions: proof of Theorem~\ref{thm:wp}(i)}
Let $f$ be a regular flux, $L=\sup|f'|$ and $u_0\in W^{3,1}(\rr)\subset\EE^2$. Let $u$ be the solution of Lemma~\ref{lem:A-mild}(a), which by Lemma~\ref{lem:A-mild}(b) is the only global mild solution, and let $T>0$.
We adapt the time-regularity argument of \cite[Ch.~6, Theorem~6.17, p.~221]{Jin} and \cite[Theorem 3.1]{JLZ18}. Here the derivatives of $u$ are limits of difference quotients, which solve equations of the form \eqref{eq:A-lin}. For $h\neq0$ let $D_hg=(g(\cdot+h)-g)/h$. Applying $D_h$ to \eqref{eq:mild} and using
\[
D_h[f(u)]=m_h\,D_hu,\qquad m_h=\int_0^1f'\big(u+\theta(u(\cdot+h)-u)\big)\,d\theta,
\]
we see that $D_hu$ solves \eqref{eq:A-lin} with $m=m_h$, the mean of $f'$ between $u$ and its translate, and $\phi_0(t)=S(t)\star D_hu_0$. The argument for $u_{xx}$ is similar, applying $D_h$ to the equation of $u_x$.
Since $u\in C([0,T+1];\BUC(\rr))$, $m_h\to f'(u)$ uniformly. Inequality \eqref{eq:A-gronwall} shows that $u$ is differentiable in $x$ and that $u_x$ solves the limit equation. 
We only write the limit equations.

\emph{Step 1. Spatial regularity.} The limits $u_x$ and $u_{xx}$ belong to $C([0,T];\EE)$ and solve \eqref{eq:A-lin} with $m=f'(u)$ and, respectively,
\[
\phi_0(t)=S(t)\star u_0',\qquad \phi_0(t)=S(t)\star u_0''-\int_0^t\partial_xP(t-s)\star [f''(u)u_x^2](s)\,ds.
\]
This is (\ref{R1}). Moreover, $A(u)\in C([0,T];\EE)$, and by Lemma~\ref{lem:A-fourier},
\begin{equation}\label{eq:A-volterra}
u(t)=u_0+\big(g_\alpha\ast A(u)\big)(t),\qquad 0\le t\le T ,
\end{equation}
pointwise, since both sides are continuous in $x$.

\emph{Step 2. Time regularity.} Since $\partial_t\widehat{S(t)}=-\xi^2\widehat{P(t)}$ and $\partial_xP(\sigma)\star f(u)=P(\sigma)\star f(u)_x$, the quotients $w_\tau=(u(\cdot+\tau)-u)/\tau$, $0<\tau\le1$, solve \eqref{eq:A-lin} with $m(s)=\int_0^1f'\big(u(s)+\theta(u(s+\tau)-u(s))\big)\,d\theta$ and
\[
\phi_0(t)=\frac1\tau\int_t^{t+\tau}P(\sigma)\star \big[u_0''-f(u)_x(t+\tau-\sigma)\big]\,d\sigma .
\]
By (\ref{K1})--(\ref{K2}), and since $g_\alpha$ is decreasing, $\|\phi_0(t)\|_\EE\le Cg_\alpha(t)$, and $\phi_0(t)\to P(t)\star A(u_0)$ in $\EE$ for $t>0$. Hence \eqref{eq:A-gronwall} with $\theta=1-\alpha$ gives $\|w_\tau(t)\|_\EE\le Ct^{\alpha-1}$ uniformly in $\tau$ and, with dominated convergence, shows that $w_\tau(t)$ converges in $\EE$, for every $t>0$, to the solution $v$ of \eqref{eq:A-lin} with
\[
m=f'(u),\qquad \phi_0(t)=P(t)\star A(u_0),
\]
and that $\|v(t)\|_\EE\le Ct^{\alpha-1}$. This gives $u\in C^1((0,T];\EE)$ with $\partial_tu=v$, which is (\ref{R2}).

By \eqref{eq:A-volterra}, $1\ast v=u-u_0=g_\alpha\ast A(u)$. Convolving with $g_{1-\alpha}$ and differentiating gives $g_{1-\alpha}\ast v=A(u)$, which together with \eqref{eq:A-volterra} is (\ref{R3}).

\subsection{Data in \texorpdfstring{$L^1\cap L^\infty$}{L1 and Linfty}: proof of Theorem~\ref{thm:wp}(ii)--(iii)}
The approximation data introduced in Remark~\ref{stantard.data}   converge to the datum in $L^1(\rr)$, and Young's inequality and the truncation give
\[
\begin{gathered}
\|u_{0,n}\|_{L^p(\rr)}\le\|u_0\|_{L^p(\rr)},\qquad -\|(u_0)_-\|_{L^\infty(\rr)}\le u_{0,n}\le\|(u_0)_+\|_{L^\infty(\rr)},\\
\|u_{0,n}-v_{0,n}\|_{L^1(\rr)}\le\|u_0-v_0\|_{L^1(\rr)},
\end{gathered}
\]
where $v_{0,n}$ are the approximations for $v_0$.

Let $R$, $(f_n)$ and $(v_{0,n})$ be as in (iii). By (i) and Lemma~\ref{lem:wp.basic}, $v_n$ is a regular and a mild solution with flux $f_n$, and
\[
|v_n|\le R,\qquad \|v_n(t)\|_{L^1(\rr)}\le\|v_{0,n}\|_{L^1(\rr)} .
\]
Lemma~\ref{lem:A-mild}(b), with $L\le L_f(R)$ and $\delta\le\delta_n(R)+\delta_m(R)$, shows that $(v_n)$ is a Cauchy sequence in $C([0,T];L^1(\rr))$ for every $T$.
By \eqref{eq:adm}, its limit $u$ is a global mild solution with flux $f$, and $|u|\le R$. By Lemma~\ref{lem:A-mild}(b) with $\delta=0$ it is the only one, so it does not depend on the choice of $R$ and of the approximations.
The convergence in $L^p(\rr)$, $p<\infty$, follows from
\[
\|v_n-u\|_{L^p(\rr)}^p\le(2R)^{p-1}\|v_n-u\|_{L^1(\rr)} .
\]
The bounds and the contraction in (ii) follow from Lemma~\ref{lem:wp.basic}, by $L^1$ convergence and the a.e.\ convergence of a subsequence. Here the lemma is applied, with the same flux $f_n$, to the regular solutions with the standard data of $u_0$ and $v_0$, and with $R$ not smaller than the $L^\infty$ norms of $u_0$ and $v_0$.

\section{Proofs of the results of Section~\ref{sec:wp.general}}\label{app:general}

In this appendix we prove Theorem~\ref{thm:wp.general}. The proof is needed for the nonlinear terms (for $F=Q=0$, see \cite[Theorem 3.1]{Zacher09}), and it constructs the approximations $u_n$ of Theorem~\ref{thm:wp.general}.
We truncate the nonlinearities (Step 1), solve the stationary problems (Step 2) and the problems with the Yosida kernels (Step 3), pass to the limit (Step 4), and prove uniqueness (Step 5). Step 6 identifies the mild solution of \eqref{convection}.

\emph{Step 1. Reductions.} If $u_0=0$, then $u=0$ is a weak solution, and it is the only one by the uniqueness proof below. Let $u_0\neq0$, $R_0=\|u_0\|_{L^\infty(\rr)}$, and let $[\,\cdot\,]_{R_0}$ be the odd $C^1$ truncation given, with $r=|\sigma|-R_0$, by
\[
[\sigma]_{R_0}=\begin{cases}
\sigma, & |\sigma|\le R_0,\\
\sgn(\sigma)\big(R_0+r-\tfrac{r^2}4\big), & R_0\le|\sigma|\le R_0+2,\\
(R_0+1)\sgn(\sigma), & |\sigma|\ge R_0+2,
\end{cases}
\]
so that
\[
0\le\frac{d}{d\sigma}[\sigma]_{R_0}\le1,\qquad \big|[\sigma]_{R_0}\big|\le R_0+1 .
\]
Since $[\sigma]_{R_0}$ has the sign of $\sigma$ and $|[\sigma]_{R_0}|\le|\sigma|$, the functions
\[
F(t,x,[\sigma]_{R_0}),\qquad Q(t,x,[\sigma]_{R_0})
\]
satisfy (\ref{hyp:F}) and (\ref{hyp:Q}) with $L_R$ replaced by $L=L_{R_0+1}$ for every $R$. A weak solution of the modified problem with $|u|\le R_0$ is a weak solution of \eqref{eq:general}.

\emph{Step 2. The stationary problems.} Let $\mathcal A(t)$ be as in Section~\ref{sec:wp.general}, with the modified $F$ and $Q$, and let $n_0=L+L^2/(2\lambda)$. For $v,\tilde v\in H^1(\rr)$ and $w=v-\tilde v$,
\begin{equation}\label{eq:monot}
\begin{split}
\langle\mathcal A(t)v-\mathcal A(t)\tilde v,w\rangle&\ge\lambda\|w_x\|^2_{L^2(\rr)}-L\|w\|_{L^2(\rr)}\|w_x\|_{L^2(\rr)}-L\|w\|^2_{L^2(\rr)}\\
&\ge\frac\lambda2\|w_x\|^2_{L^2(\rr)}-n_0\|w\|^2_{L^2(\rr)}.
\end{split}
\end{equation}
Hence, for $n>n_0$, the operator $n+\mathcal A(t)$ is strongly monotone and Lipschitz from $H^1(\rr)$ to $H^{-1}(\rr)$, and for $g\in L^2(\rr)$ the equation
\[
nv+\mathcal A(t)v=g
\]
has a unique solution $v=J(t)g\in H^1(\rr)$. It is the limit of the iteration
\[
v\mapsto v-\eps\mathcal R\big((n+\mathcal A(t))v-g\big),
\]
with $\mathcal R:H^{-1}(\rr)\to H^1(\rr)$ the Riesz isomorphism and $\eps>0$ small, and
\[
\|J(t)g-J(t)\tilde g\|_{L^2(\rr)}\le(n-n_0)^{-1}\|g-\tilde g\|_{L^2(\rr)} .
\]
Moreover, for $g\in L^2(\rr)$ and a constant $\ell>0$,
\[
\text{(S1)}\quad n\|J(t)g\|_{L^1(\rr)}\le\|g\|_{L^1(\rr)}\ \text{ if }g\in L^1(\rr),\qquad\qquad\text{(S2)}\quad |J(t)g|\le\ell\ \text{ if }|g|\le n\ell.
\]
To prove them, let $v=J(t)g$. The function $\Psi=a\,v_x-F(t,x,v)$ satisfies
\[
\Psi_x=nv-Q(t,x,v)-g\in L^2(\rr),
\]
so $\Psi$ and $v$ are continuous and vanish at $\pm\infty$.
The function $\Psi$ has the right sign at the ends of every component $(x_1,x_2)$ of $\{v>0\}$:
\[
\Psi(x_1)\ge0\ge\Psi(x_2).
\]
Indeed, if $\Psi(x_2)>0$, then $a\,v_x=\Psi+F(t,x,v)>0$ near $x_2$, which contradicts $v>0=v(x_2)$, and similarly at $x_1$.
We integrate $\Psi_x$ over the component (over bounded subintervals, if it is unbounded), where $nv-Q(t,x,v)\ge nv\ge0$ by (\ref{hyp:Q}), and we sum over the components. This gives
\[
n\int_\rr v_+\,dx\le\int_{\{v>0\}}g\,dx,
\]
and with $\{v<0\}$ we obtain (S1).
For (S2), the same argument on the components of $\{v>\ell\}$, which are bounded since $v$ vanishes at infinity and $\ell>0$, with $\Psi+F(t,x,\ell)$ in place of $\Psi$, shows that the integral of
\[
nv-g-Q(t,x,v)+\partial_xF(t,x,\ell)
\]
over them is nonpositive, although this function is positive there, by (\ref{hyp:Q}) and (\ref{hyp:F}). So $v\le\ell$, and likewise $v\ge-\ell$.

\emph{Step 3. The regularized problems.} Fix $T>0$ and $n>n_0$. With $k_n$ as in Lemma~\ref{lem:yosida}, consider
\begin{equation}\label{eq:regfix}
u_n(t)=J(t)\Big[k_n(t)u_0+\int_0^tu_n(s)\big(-k_n'(t-s)\big)\,ds\Big]\qquad\text{for a.e.\ }t\in(0,T).
\end{equation}
Since $k_n(t)+\int_0^t(-k_n'(t-s))\,ds=k_n(0)=n$, identity \eqref{eq:Rk} shows that \eqref{eq:regfix} is \eqref{eq:reg}.
Let
\[
\mathcal X=\big\{u\in L^\infty(0,T;L^2(\rr)):\ |u|\le R_0,\ \|u(t)\|_{L^1(\rr)}\le M\ \text{for a.e. }t\big\},
\]
with the distance
\[
d(u,\tilde u)=\operatorname*{ess\,sup}_{t\in(0,T)}e^{-bt}\|u(t)-\tilde u(t)\|_{L^2(\rr)} .
\]
By (S1)--(S2) with $\ell=R_0$, the right-hand side of \eqref{eq:regfix} maps $\mathcal X$ into itself.
 By the $L^2$ bound for $J(t)$, it is a contraction for $b$ large, since
\[
\int_0^Te^{-bs}\big(-k_n'(s)\big)\,ds\to0\qquad\text{as }b\to\infty .
\]
Its fixed point $u_n\in\mathcal X$ satisfies \eqref{eq:reg}, since $\mathcal A(t)u_n$ is the same for the original and the modified nonlinearities when $|u_n|\le R_0$.

We test \eqref{eq:reg} with $u_n(t)$. By \eqref{eq:Rk}, the second inequality in \eqref{eq:elem} and Lemma~\ref{lem:coercive} with $p=2$,
\[
\tfrac12\partial_t\big(k_n\ast(\|u_n\|^2_{L^2(\rr)}-\|u_0\|^2_{L^2(\rr)})\big)+\lambda\|\partial_xu_n\|_{L^2(\rr)}^2\le0 .
\]
Integrating over $(0,T)$ and using Lemma~\ref{lem:yosida}(a), we obtain
\begin{equation}\label{eq:H1bound}
2\lambda\int_0^T\|\partial_xu_n(t)\|_{L^2(\rr)}^2\,dt\le(1\ast k_n)(T)\|u_0\|^2_{L^2(\rr)}\le \Pi_\alpha(1\ast g_{1-\alpha})(T)\|u_0\|^2_{L^2(\rr)}.
\end{equation}

\emph{Step 4. The limit.} By \eqref{eq:H1bound}, the functions
\[
\psi_n=\partial_t\big(k_n\ast(u_n-u_0)\big)=-\mathcal A(t)u_n
\]
are bounded in $L^2(0,T;H^{-1}(\rr))$. By Lemma~\ref{lem:yosida}(b), $u_n-u_0=n^{-1}\psi_n+g_\alpha\ast\psi_n$. Hence, for $0<h<T$,
\begin{align*}
&\|u_n(\cdot+h)-u_n\|_{L^2(0,T-h;H^{-1}(\rr))}\\
&\qquad\le\Big(\frac2n+\|g_\alpha(\cdot+h)-g_\alpha\|_{L^1(0,T)}+\|g_\alpha\|_{L^1(0,h)}\Big)\sup_j\|\psi_j\|_{L^2(0,T;H^{-1}(\rr))}.
\end{align*}
The right-hand side is small for large $n$ and small $h$. For each fixed $n$ the left-hand side tends to $0$ as $h\to0$, by the continuity of translations in $L^2(0,T;H^{-1}(\rr))$.
For every bounded interval $I$, the space $H^1(I)$ is compactly embedded in $L^2(I)\subset H^{-1}(I)=(H^1_0(I))'$. With \eqref{eq:H1bound}, \cite[Theorem 5]{Simon87} shows that $(u_n)$ is relatively compact in $L^2(0,T;L^2(I))$.
For a subsequence, $u_n\to u$ a.e.\ in $(0,T)\times\rr$ and weakly in $L^2(0,T;H^1(\rr))$. Then $|u|\le R_0$, and Fatou's lemma gives $\|u(t)\|_{L^1(\rr)}\le M$ for a.e.\ $t$.
We integrate $\psi_n+\mathcal A(t)u_n=0$ against a test function $\varphi$ as in Definition~\ref{def:weak}. Since $k_n\ast(u_n-u_0)$ vanishes at $t=0$,
\[
\int_0^T\!\!\int_\rr\Big(-\big(k_n\ast(u_n-u_0)\big)\partial_t\varphi+a\,\partial_xu_n\varphi_x-F(t,x,[u_n]_{R_0})\varphi_x-Q(t,x,[u_n]_{R_0})\varphi\Big)\,dx\,dt=0 .
\]
Since $k_n\to g_{1-\alpha}$ in $L^1(0,T)$, the first term converges weakly in $L^2((0,T)\times\rr)$ to the corresponding term for $u$. The nonlinear terms are bounded in $L^2((0,T)\times\rr)$ and converge a.e., hence weakly.
So $u$ is a weak solution of the modified problem on $(0,T)$, and since $|u|\le R_0$, of \eqref{eq:general}. By uniqueness, these solutions do not depend on $T$, and every subsequence of $(u_n)$ has a subsequence that converges to $u$ a.e.

\emph{Step 5. Uniqueness.} Let $u$ and $\tilde u$ be weak solutions, modify $F$ and $Q$ as in Step 1 with $R_0$ larger than their $L^\infty$ norms, and let $w=u-\tilde u$. By Definition~\ref{def:weak}, $g_{1-\alpha}\ast w\in H^1(0,T;H^{-1}(\rr))$ vanishes at $t=0$, and
\[
\big\langle\partial_t(g_{1-\alpha}\ast w),w\big\rangle+\big\langle\mathcal A(t)u-\mathcal A(t)\tilde u,w\big\rangle=0\qquad\text{for a.e.\ }t\in(0,T).
\]
The uniqueness part of the proof of \cite[Theorem 3.1]{Zacher09}, with $H^1(\rr)$ and $L^2(\rr)$ as the spaces $V$ and $H$ there, uses the equation only in this form, together with a G\aa rding inequality for the second term. Here, by \eqref{eq:monot},
\[
\langle\mathcal A(t)u-\mathcal A(t)\tilde u,w\rangle\ge\frac\lambda2\|w\|^2_{H^1(\rr)}-\Big(n_0+\frac\lambda2\Big)\|w\|^2_{L^2(\rr)} .
\]
That argument, which uses the Yosida approximations of Lemma~\ref{lem:yosida} and a Gronwall lemma, therefore gives $w=0$.

\emph{Step 6. The mild solution of \eqref{convection}.} The regular solutions $v_n$ of Theorem~\ref{thm:wp}(iii), with the standard data, satisfy Definition~\ref{def:weak} with $a\equiv1$, $Q=0$, the fluxes $f_n$ and these data, by (\ref{R3}).
Since the convection drops out, \eqref{acv2} and the second inequality in \eqref{eq:elem} give
\[
2\int_0^T\|\partial_xv_n(t)\|_{L^2(\rr)}^2\,dt\le(1\ast g_{1-\alpha})(T)\|u_0\|_{L^2(\rr)}^2,
\]
and the identity of Definition~\ref{def:weak} passes to the limit. So the mild solution is a weak solution, and by uniqueness it is the weak solution.

\section*{Acknowledgments}
This work started in 2015, when the author visited Mahamadi Warma at the University of Puerto Rico and discovered fractional-in-time derivatives. The dyadic energy method for the linear problem appeared some years later. Over the years, despite many discussions with colleagues, among them Timoth\'ee Crin-Barat,  Grzegorz Karch, Julio Rossi, Fernando Quir\'os and Enrico Valdinoci, we did not succeed in applying it to nonlinear problems. It was in the fall of 2026 that, with the help of Anthropic Claude Opus, the author found the simpler argument presented in the introduction. The first energy method for nonlinear problems used small initial data to obtain the correct decay, an assumption suggested at some point by Timoth\'ee Crin-Barat, inspired by his previous works. The initial result was then improved to remove this assumption.

\section*{AI Disclosure}
The author used Anthropic Claude Opus  to help develop and check the proofs and to prepare the manuscript. After using this tool, the author reviewed and edited the content as needed. The author takes full responsibility for the content of the publication.

\end{document}